\documentclass[aap, seceqn,dvips]{arximspdf}

\doi{10.1214/07-AAP459}
\volume{18}
\issue{2}
\pubyear{2008}
\firstpage{427}
\lastpage{469}

\makeatletter
\renewcommand{\theequation}{\arabic{section}.\arabic{equation}}
\newproclaim{notation}{Notation}
\newtheorem{lema}{Lemma}
\newtheorem{teorema}{Theorem}
\newproclaim{cas}{Case}
\newtheorem{proposition}{Proposition}
\newproclaim{remark}{Remarks}
\newproclaim{examples}{Examples}
\newproclaim{step}{Step}
\newcommand{\vid}{\nonumber \\[-8pt]\\[-8pt]}
\makeatother

\begin{document}
\begin{frontmatter}

\title{Functional quantization rate and mean regularity of
processes with an application\break to L\'evy processes}
\runtitle{Functional quantization and regularity of processes}

\begin{aug}
\author[A]{\fnms{Harald} \snm{Luschgy}\ead[label=e1]{luschgy@uni-trier.de}\corref{}}
and
\author[B]{\fnms{Gilles} \snm{Pag\`es}\ead[label=e2]{gpa@ccr.jussieu.fr}}
\runauthor{H. Luschgy and G. Pag\`es}
\affiliation{Universit\"at Trier and Universit\'e Paris 6}
\address[A]{
Universit\"at Trier\\
FB IV-Mathematik\\
D-54286 Trier\\
BR Deutschland\\
\printead{e1}} %adresu isvedimo komanda gale!
\address[B]{
Laboratoire de Probabilit\'es\\
\quad et Mod\`eles al\'eatoires\\
UMR~7599\\
Universit\'e Paris 6\\
case 188\\
4 pl. Jussieu\\
F-75252 Paris Cedex 5\\
France\\
\printead{e2}}
\end{aug}

% HISTORY:
\received{\smonth{2} \syear{2007}}
\revised{\smonth{2} \syear{2007}}

% ABSTRACT
\begin{abstract}
We investigate the connections between the mean pathwise regularity
of stochastic processes  and their $L^r(\mathbb{P}$)-functional quantization
rates  as random variables taking
values in some $L^p([0,T],dt)$-spaces ($0<p\le r$). Our main tool is
the Haar basis. We then emphasize that the derived functional
quantization rate may be optimal
(e.g.,  for Brownian motion or symmetric stable
processes) so that the rate is optimal as a universal upper bound. As a
first application, we establish the $O((\log N)^{-1/2})$ upper
bound for general It\^o
processes which include multidimensional diffusions. Then, we focus
on the specific family of L\'evy processes for which we derive a
general quantization rate
based on the regular variation properties of its L\'evy measure at~0.
The case of compound Poisson processes, which appear as
degenerate in the former approach, is studied specifically: we
observe some rates which are between the
finite-dimensional and infinite-dimensional ``usual'' rates.
\end{abstract}

% KEYWORDS
\begin{keyword}[class=AMS]
\kwd{60E99}
\kwd{60G51}
\kwd{60G15}
\kwd{60G52}
\kwd{60J60}.
\end{keyword}
\begin{keyword}
\kwd{Functional  quantization}
\kwd{Gaussian process}
\kwd{Haar basis}
\kwd{L\'evy process}
\kwd{Poisson process}.
\end{keyword}
\end{frontmatter}
%
%s1 ###
\section{Introduction}
In this paper, we investigate the connection
between the functional
$L^r(\mathbb{P})$-quantization rate for a process $X=(X_t)_{t\in[0,T]}$ and
the  $L^r(\mathbb{P})$-mean pathwise regularity
of the mapping $t\mapsto X_t$ from
$[0,T]\to L^r(\mathbb{P})$ in an abstract setting by means of a constructive
approach (we mean that all the rates
are established using some explicit sequences of quantizers).

First, let us briefly recall what functional  quantization is and how
it was introduced. Let $(E,\|\cdot\|)$
denote a finite-dimensional ($E=\mathbb{R}$ or
$\mathbb{R}^d$)  or infinite-dimensional ($E= L^p([0,T],dt)$, $1 \leq
p<\infty$, ${\mathcal{C}}([0,T])$,\,\ldots) separable Banach space (or
complete quasi-normed space like $E=
L^p([0,T],dt)$, $0<p<1$) and let
$\alpha \subset E$ be a finite subset of size
$\operatorname{card}(\alpha) \leq N$, $N\geq 1$. The \textit{ Voronoi quantization} of
an $E$-valued random vector $X\dvtx(\Omega,\mathcal{A},\mathbb{P})\to E$ \textit{with
respect to the
codebook  $\alpha$} is simply the projection of $X$ onto $\alpha$ following
the nearest neighbor rule, that is
\[
\widehat X^\alpha = \pi_\alpha(X),
\]
where
\[
\pi_\alpha = \sum_{a\in \alpha}a\mathbf{1}_{C_a(\alpha)},
\]
$(C_a(\alpha))_{a\in\alpha}$ being a Borel partition of $E$ satisfying, for every
$a\!\in\alpha$,
\[
C_a(\alpha) \subset \biggl\{u\!\in E\dvtx\|u-a\|\le
\min_{b\in\alpha\setminus\{a\}}\|u-b\|\biggr \}.
\]
Then, the $L^r$-mean quantization error ($0<r<\infty$) is defined by
\[
\|X-\widehat X^\alpha\|_{L_E^r(\mathbb{P})}= \biggl(\mathbb{E} \min_{a\in
\alpha}\|X-a\|^r\biggr)^{1/r}.
\]
This quantity is finite provided $X\!\in L_E^r(\mathbb{P})$. The  set
$\alpha$ is called an $N$-\textit{codebook} or $N$-\textit{quantizer}. It can be shown that
such random vectors
$\widehat X^\alpha$ are the best approximation of $X$ among all
$\alpha$-valued random vectors. The minimal $N$th quantization error
of $X$ is then defined by
%
%e1.1 ###
\begin{equation}
e_{N,r}(X,E) := \inf \biggl\{\biggl(\mathbb{E} \min_{a \in \alpha} \| X-a \|^r
\biggr)^{1/r} \dvtx \alpha \subset E, \operatorname{card}(\alpha) \leq N \biggr\}.
\end{equation}
When $E= L^p([0,T],dt)$ (with its usual norm or quasi-norm denoted by
$|\cdot|_{L^p_{_T}}$ from now on), an
$E$-valued random variable $X$ is  a (bimeasurable) stochastic
process $X=(X_t)_{t\in[0,T]}$
defined on the probability space $(\Omega,\mathcal{A}, \mathbb{P})$ whose
trajectories $(X_t(\omega))_{0\le t\le T}$ (almost) all belong to
$L^p([0,T],dt)$. The $L^r$-integrability
assumption then reads
\[
\mathbb{E}\biggl(\biggl(\int_0^T|X_t|^p\,dt\biggr)^{r/p}\biggr)<+\infty.
\]

It is still an open question whether
$L^r$-optimal
$N$-quantizers for Gaussian random vectors always exist in an
abstract Banach space setting  (see~\cite{GRAF3}). However, in many
situations of interest for processes, including all the\break
$L^p([0,T],dt)$-spaces, $1\le p<+\infty$, the existence of at least
one
such
$L^r$-optimal codebook has been established (provided $\mathbb{E} \|X\|^r
<+\infty$).  Note, however, that this is
not the case for the space $\mathcal{ C}([0,T])$ of continuous
functions. For more details on the existence problem for optimal
quantizers, we refer to~\cite{GRAF3}.

On the other hand, optimal $L^r$-quantizers always exist when
$E=\mathbb{R}^d$, $d\geq 1$. In this finite-dimensional
setting,   this problem is known as \textit{optimal vector
quantization} and has been
extensively investigated since the early 1950s  with some
applications to signal processing and transmission
(see~\cite{IEEE} or~\cite{GEGR}).

In $d$-dimensions, the convergence rate  of $e_{N,r}$ is given by the
so-called Zador
theorem,
%e1.2 ###
\begin{equation}\label{Zador}
\lim_N N^{1/d} e_{N,r}(X,\mathbb{R}^d) = \tilde{J}_{r,d}\biggl(\int_{\mathbb{R}^d}
g^{{d}/{(d+r)}}(\xi)\,d\xi\biggr)^{1/r +1/d},
\end{equation}
where $g$ denotes the density of the absolutely continuous part of
the distribution $\mathbb{P}_{X}$ of $X$ and $\tilde{J}_{r,d}\in(0,\infty)$ (see~\cite{GRLU}).

Since the early 2000's, much attention has been paid to the infinite-dimensional
case. This is the so-called
\textit{functional quantization} problem for  stochastic processes: the
aim is to quantize some processes
viewed as random vectors taking values in their path spaces, supposed
to be $
L^p([0,T],dt)$ spaces,
$1\le p<+\infty$. Many results have been obtained for several
families of processes with special attention having been paid to
Gaussian processes and (Brownian) diffusion processes by several
authors.  Thus, in the purely Hilbert space setting ($r=2$, $E=
L^2([0,T],dt)$), the sharp rate of quantization of the Brownian motion
$(W_t)_{t\in[0,T]}$ is given (see~(3.6) in~\cite{LUPA2}) by
%e1.3 ###
\begin{equation}\label{sharpMB}
e_{N,2}(W,L^2([0,T],dt)) \sim \frac{\sqrt{2}\,T}{\pi(\log N)^{
1/2}}.
\end{equation}
The existence of such  a \textit{sharp rate} for Brownian motion has
been extended to $L^p([0,T],dt)$ spaces for $1\le p\le \infty$
(see~\cite{Dereichetal}). Similar sharp rates (with an explicit
constant)  hold for a wide class of Gaussian processes, including the
fractional Brownian motions for which we have
\[
e_{N,2}(W^H,L^2([0,T],dt)) \sim \frac{c(H,T)}{(\log N)^{H}},
\]
where $H$ denotes the Hurst parameter of the fractional Brownian
motion $W^H$, the Ornstein--Uhlenbeck process, the Brownian
sheet, and so on, in the purely Hilbert space setting (see~\cite{LUPA2}). The
exact rate has also been established in~\cite{LUPA1} (Section~3)
for a wider class of Gaussian processes. In ~\cite{LUPA1,LUPA2},
these results are based on the (sharp or exact) asymptotic
behavior of the eigenvalues  of high order of the Karhunen--Lo\`eve
expansion of the Gaussian process. As a byproduct, this
approach provides very simple explicit sequences of rate-optimal
asymptotic quantizers (provided that the Karhunen--Lo\`eve
expansion of the process itself is accessible). Their numerical
implementation has lead to some unexpectedly promising
numerical applications in finance, especially for the pricing of
path-dependent options like  Asian options in several
popular models of asset dynamics (Black--Scholes, stochastic
volatility Heston and  SABR models, etc.). For these
aspects, we refer to~\cite{PAPR2} or~\cite{WI}. We also mention
  applications of quantization to statistical
clustering of  data  (see, e.g.,~\cite{POST}) and
some more recent developments concerning functional data
investigated in~\cite{TAKI} and~\cite{TAPEOG}.

For Gaussian processes, an important connection
with the small ball probability problem has been made
(see~\cite{DEFEMASC,GRLUPA1}). Some exact or sharp rates of
convergence for
different classes of Brownian diffusions have also recently been
proven (see~\cite{Dereich,LUPA4})
with a rate driven by $(\log N)^{-1/2}$.

The common feature shared by all these results is that there is a
one-to-one correspondence between the exponent $a$
that controls the $(L^r(\mathbb{P}),L^p(dt))$-quantization rate of these
processes in the  $\log (N)^{-a}$ scale and their mean
pathwise regularity, that is, the largest exponent $a$ that satisfies
%
%e1.4 ###
\begin{equation}\label{meanpathcont}
\forall\, s,\, t\!\in [0,T]\qquad \|X_t-X_s\|_{L^r(\mathbb{P})}\leq C_r|t-s|^a.
\end{equation}

Although such a correspondence is not really surprising given the
connection between quantization rate and small ball
probabilities in the Gaussian setting, this  naturally leads to an
attempt to derive
a general abstract result that connects these two
features of a process. This is the aim of
Section~\ref{upper-bound} of this paper, in which  we show that the
mean pathwise regularity always provides a universal upper bound for
the $(L^r(\mathbb{P}), L^p(dt))$-quantization rate ($0\!<\!p\!\le\! r$). We
then retrieve
the rate obtained by more specific approaches for all the
processes mentioned above. We also extend  to general
Brownian diffusion processes and even general It\^o processes the
rate formerly obtained for specific classes of diffusions
in~\cite{Dereich,LUPA4}. We also obtain some first quantization
rates for some classes of L\'evy processes.  The main technique is
to expand a process on the simplest  wavelet basis---the Haar basis
(known to be unconditional when $p>1$)---and to use a
nonasymptotic  version of the Zador theorem (a slight
improvement  of the Pierce lemma; see~\cite{GRLU}).

At this point, the  next question is to ask conversely whether
this always provides the true quantization rate. In
this na\"{i}ve form, the answer to this question is clearly ``no'' because
equation~(\ref{meanpathcont}) only takes into account
the  mean pathwise  H\" older regularity of a process and one can
trivially build (see~\cite{LUPA1}) some processes with
smoother mean regularity (like processes with $\mathcal{C}^k$, $k\ge 1$,
trajectories). We do not extend our approach
in that direction, for the sake of simplicity, but there is no doubt
that developing  techniques similar to those used in
Section~\ref{upper-bound}, one can connect higher order mean pathwise
regularity and quantization rate, as in the H\"older
setting.   This would require an appropriate wavelet basis.
In fact, we point out  in Section~\ref{LevyLevy}, devoted to general
L\'evy processes, that the answer may be
negative---the quantization rate can be infinitely faster
than the
mean pathwise regularity---for different reasons in
connection with the dimensionality of the process: a Poisson process
is,
in some sense, an almost finite-dimensional random vector
which induces a very fast quantization rate which does not take place
in the $(\log N)^{-a}$ scale, although the mean pathwise
$L^r(\mathbb{P})$-regularity of a Poisson process is H\"older [and depends on
$r$; see, e.g.,~(\ref{Poissrle1})
and~(\ref{Poissrge1})].  Conversely, we emphasize, via on several classes
of examples, that the upper bound derived from
mean regularity   provides the true rate of quantization. This
follows from  a comparison  with the lower bound that
can be  derived from small deviation results
(see,~e.g.,~\cite{GRLUPA1} or the remark below Theorem~\ref{Upper1}
which elucidates the connection between functional quantization and
small deviation theory). Thus, we prove
that our approach yields the exact rate for a wide class of
subordinated L\'evy processes (including symmetric
$\alpha$-stable processes).

The main result of  Section~\ref{LevyLevy} is Theorem~\ref{Levyrate},
which provides a  functional  quantization rate for a general
L\'evy process $X$ having no Brownian component: this rate is
controlled
by the behavior of the L\'evy measure
$\nu$ around $0$ (e.g., the index of $X$ for a stable process). As
an example for L\'evy processes which have
infinitely many small jumps,  if the (infinite) L\'evy measure $\nu$
(is locally absolutely continuous  around
0 and)   satisfies
\[
\exists\, c>0\qquad\mathbf{1}_{\{0<|x| \leq c]}\nu(dx)\leq\frac{ C
}{|x|^{\underline \theta +1}}\mathbf{1}_{\{0<|x| \leq c]}\,dx
\]
for some $\underline \theta\!\in(0,2]$, then, for every $p,\,
r\!\in(0,\underline \theta]$ such that $0<p\leq r$ and $X_1\!\in
L^r(\mathbb{P})$,
\[
e_{N,r}(X,L^p([0,T],dt))  = O((\log N)^{-{1}/{\underline \theta}}).
\]
This makes a connection between quantization rate and   the
Blumenthal--Getoor index $\beta$
of $X$ when
$\nu$  satisfies the above upper bound with $\underline\theta
=\beta$. In fact,  a more general result is
established in Theorem~\ref{Levyrate}:  when the ``$0$-tail
function'' $\underline \nu \dvtx x\mapsto
\nu([-x,x]^c)$ has regular variation   as $x$ goes to $0$, with index
$-\underline \theta$, then   $\underline \theta=\beta$
(see~\cite{BLGE}) and we establish a close  connection between  the
quantization rate of $X$ and $\underline \nu$, $\underline
\theta$.   In many cases of interest, including $\alpha$-stable
processes and other classes of subordinated L\'evy
processes, we show that this general upper bound provides  the exact
rate of quantization; it matches the lower bound
estimates derived from the connection between quantization rate and
small deviation estimates (see, e.g.,~\cite{GRLUPA1}). When
the L\'evy process does have a Brownian component, its exact
quantization rate is $ (\log N)^{-{1}/{2}}$, like Brownian
motion [when $0<p<r<2$, $X_1\!\in L^r(\mathbb{P})$].

When the L\'evy measure is finite (then  $\underline \theta =0$), we also
establish some quantization rates for the compound Poisson processes
and show they are infinitely faster than the above ones. To
this end, we design  an explicit sequence of quantizers which can clearly be implemented
for numerical purposes.  In fact, the whole
proof is constructive, provided the L\'evy measure is ``tractable'' enough.

The paper is organized as follows. Section~\ref{upper-bound} is
devoted to the abstract connection between mean regularity and
quantization
rate of processes. Section~\ref{FirstAppli} is devoted to some initial
applications to various  families of processes.  As far as we know,
some of
these rates are new. In  several cases of interest,  these rates are
shown to be optimal. The main result is
Theorem~\ref{Upper1}. Section~\ref{LevyLevy} provides an upper bound
for the quantization rate of general L\'evy process in connection
with the behavior of
the L\'evy measure around~$0$. The main results are
Theorem~\ref{Levyrate} and Proposition~\ref{Poisson}.  In
Section~\ref{LevyB}, we provide the exact rate for a L\'evy process
having a Brownian component. Finally, in
Section~\ref{suboLevy}, we derive the exact quantization rate for
subordinated L\' evy processes.

\begin{notation*}
\begin{itemize}
\item $L^p_{T}:=L^p([0,T],dt)$ and
$|f|_{L^p_{T}}= (\int_0^T
|f(t)|^p\,dt)^{1/p}$.

\item Let $(a_n)_{n\ge 0}$ and $(b_n)_{n\ge 0}$ be two
sequences of positive real numbers.
$a_n\sim b_n$ means $ a_n=b_n +o(b_n)$ and $a_n\approx b_n$ means
$a_n =O(b_n)$ and
$b_n=O(a_n)$.

\item $[x]$ denotes the integral part of the real
number $x$ and $x_{+}=\max(x,0)$ its positive part.

\item $\log_m(x)$ denotes the $m$-times iterated
logarithm function.

\item $\|Y\|_{r}:= \|Y\|_{L^r(\mathbb{P})}$ for any random
variable $Y$ defined on a probability space $(\Omega,\mathcal{A}, \mathbb{P})$.

\item Throughout the paper, the letter $C$ (possibly
with subscripts) will denote a positive real constant that may vary
from line to line.

\item For a c\`adl\`ag continuous-time process
$X=(X_t)_{t\ge 0}$, $X_{t-}$ will denote its left limit and $\Delta
X_t :=X_t-X_{t-}$ its jump at time $t$.
\end{itemize}
\end{notation*}

%s2 ###
\section{Mean pathwise regularity and quantization error rate: an
upper bound}\label{upper-bound}
In this section, we   derive in full generality an upper bound for the
$(L^r(\mathbb{P}),L^p_{T})$-quantization error $e_{N,r}(X,L^p_{T})$ based
on the path regularity
of the mapping $t\mapsto X_t$ from $[0,T]$ to $L^\rho(\mathbb{P})$. The main
result of this section is
Theorem~\ref{Upper1} below. We will then  illustrate via several
examples  that this rate
may be optimal or not.

As a first step, we will reformulate the so-called Pierce lemma
(see~\cite{GRLU}, page 82), which is
the main step of  the proof of Zador's Theorem for unbounded random
variables. Note that the proof of its original formulation (see
below)
relies   on random quantization.
\begin{lema}[\textup{(Extended Pierce Lemma)}]\label{UPierce}
Let $r,\delta>0$. There exists a real
constant $C_{r,\delta}$
such that, for every random
variable
$X\dvtx(\Omega,\mathcal{A})\to (\mathbb{R},\mathbb{B}(\mathbb{R}))$,
\[
\forall\,N\ge 1
\qquad e_{N,r}(X,\mathbb{R})= \inf_{{\operatorname{card}}(\alpha)\leq N}
\|X-\widehat
X^\alpha\|_{r} \leq
C_{r,\delta}\|X\|_{{r+\delta}}\,N^{-1}.
\]
\end{lema}
\begin{pf}
It follows from the original
Pierce lemma that there exists
a universal real constant $C^0_{r,\delta}>0$
and  an integer
$N_{r,\delta}\ge 1$ such  that, for any random
variable
$X\dvtx(\Omega,\mathcal{A})\to (\mathbb{R},\mathbb{B}(\mathbb{R}))$,
\[
\forall\,N\geq
N_{r,\delta} \qquad \inf_{{\operatorname{card}}(\alpha)\leq N} \mathbb{E}|X-\widehat
X^\alpha|^r \le
C^0_{r,\delta}(1+\mathbb{E}|X|^{r+\delta})\,N^{-r}.
\]

Using the scaling property of quantization,  for every
$\lambda>0$,
\[
\|X-\widehat X^{\alpha}\|_{r}=
\frac{1}{\lambda}\|(\lambda
X)-\widehat{\lambda
X}^{\lambda\alpha}\|_{r},
\]
where
$\lambda
\alpha =\{\lambda a,\; a\!\in \alpha\}$, one derives  from the
Pierce lemma, by considering $X /\|X\|_{r+\delta}$ and setting
$\lambda
:=1/\|X\|_{{r+\delta}}$, that
\[
\forall\,N\ge
N_{r,\delta} \qquad \inf_{{\operatorname{card}}(\alpha)\le N} \|X-\widehat
X^\alpha\|_{r} \le
(2C^0_{r,\delta})^{
1/r}\|X\|_{{r+\delta}}\,N^{-1}.
\]
Now, for
every $N\!\in\{1,\ldots,N_{r,\delta}-1\}$, setting $\alpha:=\{0\}$
yields
\[
\inf_{\operatorname{card}(\alpha)\le N} \|X-\widehat X^\alpha\|_{r}
\le \|X\|_{r} \le N_{r,\delta}
\|X\|_{{r+\delta}}\,N^{-1}.
\]
Combining the last two inequalities
and setting $C_{r,\delta}\!=\!\max(
(2C^0_{r,\delta})^{
1/r},
N_{r,\delta})$ completes the proof.
\end{pf}

Let $(e_n)_{n\ge 0}$ denote the Haar basis,  defined as the
restrictions to $[0,T]$ of the following functions:
\begin{eqnarray*}
e_0&:=& T^{-1/2}\mathbf{1}_{[0,T]},\qquad e_1:=T^{-
1/2}\bigl(\mathbf{1}_{[0,T/2)}-\mathbf{1}_{[T/2,T]}\bigr),\\
e_{2^n+k}&:=& 2^{
n/2}e_1(2^n\cdot-kT),\qquad n\ge0, \, k\!\in\{0,\ldots,2^{n}-1\}.
\end{eqnarray*}

With this normalization, it makes up an orthonormal basis of the
Hilbert space $(L^2_{T}, (\cdot|\cdot))$, where
$(f|g)=\int_0^Tfg(t)\,dt$ and a (monotone) Schauder basis of
$L^p_{T}$, $p\!\in[1,+\infty)$, that is,
$(f|e_{0})e_{0}+\sum_{n\ge 0}\sum_{0\le k\le
2^n-1}(f|e_{2^n+k})e_{2^n+k}$,
converges to $f$ in
$L^p_{T}$   for every $f\!\in L^p_{T}$ (see~\cite{SI}).
Furthermore, it clearly satisfies, for every $f\!\in
L^1_{T}$ and  every
$p>0$,
%
%e2.1 ###
\begin{eqnarray}\label{incondHaar}
\qquad&&\forall\,n\ge 0
\vid
&&\qquad \int_0^T
\Bigg|\sum_{k=0}^{2^n-1}(f|e_{2^n+k})e_{2^n+k}(t)\Bigg|^p\,dt
= 2^{n(
p/2-1)}T^{1-p/2}\sum_{k=0}^{2^n-1}|(f|e_{2^n+k})|^p.\nonumber
\end{eqnarray}

The second key to establish a general connection between quantization
rate and mean pathwise
regularity is the following standard property of the Haar basis: for every
$f\!\in L^1_{T}$,
\begin{eqnarray}\label{1.2}
\qquad&&(f|e_{2^n+k})\nonumber\\
&&\qquad=2^{n/2}T^{-
1/2}\biggl(\int_{kT2^{-n}}^{(2k+1)T2^{-(n+1)}}f(u)\,du-\int_{(2k+1)T2^{-(n+1)}}
^{(k+1)T2^{-n}}f(u)\,du\biggr)
\vid
&&\qquad=2^{n/2}T^{-
1/2}\int_{0}^{T2^{-(n+1)}}\bigl(f(kT2^{-n}+u)\nonumber\\
&&\hspace{104pt}\qquad{}-f\bigl((2k+1)T2^{-(n+1)}+u\bigr)\bigr)\,du.\nonumber
\end{eqnarray}

Let $(X_t)_{t\in[0,T]}$ be a  bimeasurable process
defined on a probability space $(\Omega,\mathcal{A}, \mathbb{P})$ with
$\mathbb{P}$-almost all paths lying
in $L^1_{T}$  such that
$X_t\!\in L^{\rho}(\mathbb{P})$ for every  $t\!\in[0,T]$ for some positive
real exponent $\rho>0$. When $\rho\!\in  (0,1)$, we assume that $X$
has c\`adl\`ag  paths (right-continuous,
left-limited)  to ensure the measurability of the
supremum  in assumption~(\ref{assumption1}) below.

We make the following  $\varphi$-Lipschitz assumption on the map
$t\mapsto X_t$ from $[0,T]$
into
$L^{\rho}(\mathbb{P})$:  there is a nondecreasing function
$\varphi\!:\!\mathbb{R}_+\to [0,+\infty]$, continuous at $0$ with $\varphi
(0)\!=\!0$,  such that
%
%e2.2 ###
\begin{eqnarray}\label{assumption1}
( L_{\varphi,\rho})\equiv
\cases{
\mbox{(i)}\quad \forall\,s,\,t\!\in[0,T],\cr
\qquad\mathbb{E}\,|X_t-X_s|^\rho
\le (\varphi(|t-s|))^{\rho},
 \qquad \mbox{if $\rho\ge 1,$} \cr
\mbox{(ii)}\hspace{9.5pt} \forall\,t\in[0,T],\,\forall\,h\!\in(0,T],\cr
\qquad\mathbb{E}\biggl(\displaystyle\sup_{t\le s\le (t+h)\wedge  T}|X_s-X_t|^\rho\biggr) \le
(\varphi(h))^{\rho},
\cr
\hspace{178pt}\mbox{if  $0<\rho< 1.$}
}
\end{eqnarray}
[One may assume, without loss of generality, that $\varphi$ is always
finite, but that (i)~and (ii) are only true for $|t-s|$ or $h$
small enough, resp.] Note that this
assumption implies that
$\mathbb{E} ( |X|^\rho_{L^\rho_{T}})<+\infty$ so that, in particular,
$\mathbb{P}(d\omega)$-a.s.,
$t\mapsto \!X_t(\omega)$  lies in   $L^\rho_{T}$ (which, in turn,
implies that the paths lie in  $L^1_{T}$ if $\rho\ge 1$).

We make    a  \textit{regularly varying assumption} on $\varphi$ at
$0$ with index
$b\ge 0$, that is,  for every  $t>0$,
%e2.3 ###
\begin{equation}\label{RV}
\lim_{x\to 0} \frac{\varphi(tx)}{\varphi(x)}= t^b.
\end{equation}
In accordance with the  literature (see~\cite{BIGOTE}), this means
that $x\mapsto \varphi(1/x)$
is regularly varying at infinity with index $-b$ (which is a more
usual notion in that field). When $b=0$, $\varphi$ is said to be \textit{slowly varying}
at $0$.

Let $r,\,p\!\in (0,\rho)$. Our aim is to evaluate the $L^r(\mathbb{P})$-quantization
rate of  the process $X$, viewed as an
$L^p_{T}$-valued  random variable induced  by   the ``Haar product
quantizations'' of $X$ defined by
%
%e2.4 ###
\begin{equation}
\widehat X = \widehat{\xi}^{N_0}_0e_0 + \sum_{n\ge 0}
\sum_{k=0}^{2^n-1}\widehat{\xi}^{N_{2^n+k}}_{2^n+k}e_{2^n+k},
\end{equation}
where $\xi_{k} :=(X|e_k)\!\in L^\rho(\mathbb{P})$, $k\ge 0$, and
where $\widehat{\xi}^{N}$ denotes an $N$-quantization
($N\ge 1$) of the (real-valued) random variable $\xi$, that is, a
quantization of $\xi$ by a codebook $\alpha^N$
having $N$ elements. A quantization taking finitely many values, we
set $N_{2^n+k}=1$ and
$\widehat{\xi}^{N_{2^n+k}}_{2^n+k}=0$ for large enough $n$ (which may
be a nonoptimal~$1$-quantizer for
$\xi^{N_{2^n+k}}_{2^n+k}$).

We will see that this  local behavior of $\varphi$ at $0$ induces an
upper bound for the functional
quantization error rate of $X$ (regardless of the values of
$r$ and $p$, except for constants).
\begin{teorema}\label{Upper1}
Let  $X=(X_t)_{t\in[0,T]}$  be  a
(bimeasurable) process defined on a probability space $(\Omega,
\mathcal{A},\mathbb{P})$ such that $X_t\!\in L^\rho(\mathbb{P})$  for an exponent $\rho>0$.
Assume that $X$ satisfies~\textup{(\ref{assumption1})} [the
$\varphi$-Lipschitz assumption $(L_{\varphi,\rho})$]  for this
exponent $\rho$, where $\varphi$ is regularly varying \textup{[}in the sense
of~\textup{(\ref{RV})}\textup{]} with index $b\ge 0$ at
$0$ \textup{[}then $|X|_{L^\rho_{T}}\!\in L^1(\mathbb{P})$\textup{]}. Then
\[
\forall\, r, p\in(0,\rho)\qquad e_{N,r}(X,L^p_{T}) \le C_{r,p}
\cases{
\varphi(1/\log N), & \quad if $b>0$,\cr
\psi(1/\log N), & \quad if $b=0$,
}
\]
with $\psi(x) =(\int_0^x (\varphi(\xi))^{r\wedge 1}\,d\xi/ \xi)^{1/(r\wedge 1)}$, assuming,
moreover, that\break $\int_{0}^1
(\varphi(\xi))^{r\wedge 1}\,d\xi/ \xi <+\infty$  if $b=0$.  In
particular, if
$\varphi(u)=cu^b$,
$b>0$, then
%
%e2.5 ###
\begin{equation}\label{puissance}
e_{N,r}(X,L^p_{T}) = O((\log N)^{-b}).
\end{equation}
\end{teorema}

\begin{pf}
Using the two  obvious inequalities
\[
|f|_{L^p_{T}} \le  T^{1/p-1/p'}|f|_{L^{p'}_{T}},\qquad p\le p',
\]
for  every  Borel function $f \dvtx [0,T]\to \mathbb{R}$ and
\[
\|Z\|_r\le \|Z\|_{r'}, \qquad r\le r',
\]
for  every  random variable $Z\dvtx \Omega\to \mathbb{R}$, we may assume, without
loss of generality, that either
\[
1\le p=r <\rho\quad \mbox{or}\quad 0<p=r<\rho\le 1.
\]
\begin{cas}[($1\le p=r <\rho$)]
Let $N\ge 1$ be
a fixed integer. We consider a Haar product
quantization
$\widehat
X$ of $X$
with a (product) codebook having at most $N$ elements, that is,
such that $N_0\times \prod_{n,k} N_{2^n+k}\le N$.  Its
characteristics will be specified below.  Then,
using~(\ref{incondHaar}), that is,
\begin{eqnarray*}
|X-\widehat X|_{L^r_{T}}
&\le& T^{1/r-1/2}|\xi_0-\widehat \xi_0^{N_0}|+
\sum_{n\ge 0}
\Bigg|\sum_{k=0}^{2^n-1}(\xi_{2^n+k}-\widehat
\xi_{2^n+k}^{N_{2^n+k}})e_{2^n+k}\Bigg|_{L^r_{T}}\\
&=&T^{1/r-1/2}|\xi_0-\widehat \xi_0^{N_0}|\\
&&{}+ T^{1/r-
1/2}\sum_{n\ge 0}2^{n(1/2 -
1/r)}
\Biggr(\sum_{k=0}^{2^n-1}|\xi_{2^n+k}-\widehat
\xi_{2^n+k}^{N_{2^n+k}}|^r\Biggr)^{1/r}
\end{eqnarray*}
so that, both $\|\cdot\|_{r}$ and $\|\cdot\|_{1}$ being
norms,
%
%e2.11 ###
%e2.10 ###
%e2.9 ###
%e2.8 ###
%e2.7 ###
%e2.6 ###
\begin{eqnarray}\label{Fonda1}
&&\bigl\||X-\widehat X|_{L^r_{T}}\bigr\|_{r}\nonumber\\
&&\qquad\le T^{ 1/r-1/2} \bigl\||\xi_0-\widehat \xi_0^{N_0}|
\bigl\|_{r}\nonumber\\
&&\qquad\quad{}+   T^{
1/r-1/2}\sum_{n\ge 0}2^{n(1/2 -
1/r)}\Bigg\|
\Biggl(\sum_{k=0}^{2^n-1}|\xi_{2^n+k}-\widehat
\xi_{2^n+k}^{N_{2^n+k}}|^r\Biggr)^{
1/r}\Bigg\|_{r}\nonumber\\
&&\qquad=T^{1/r-1/2} \| \xi_0-\widehat \xi_0^{N_0}\|_{r}\nonumber\\
&&\qquad\quad{}+T^{1/r-1/2}\sum_{n\ge 0}2^{n( 1/2 -
1/r)} \Bigg\|
\sum_{k=0}^{2^n-1}|\xi_{2^n+k}-\widehat
\xi_{2^n+k}^{N_{2^n+k}}|^r
\Bigg\|_{1}^{1/r}\nonumber\\
&&\qquad\le T^{1/r-1/2}\|
\xi_0-\widehat \xi_0^{N_0} \|_{r}\\
&&\qquad\quad{}+   T^{ 1/r-
1/2}\sum_{n\ge 0}2^{n( 1/2 -
1/r)}\biggl( 2^n \max_{0\le k\le
2^n-1}\bigl\|
|\xi_{2^n+k}-\widehat
\xi_{2^n+k}^{N_{2^n+k}}|^r\big\|_{1}\biggr)^{
1/r}\nonumber\\
&&\qquad= T^{1/r-1/2}\|\xi_0-\widehat \xi_0^{N_0}
\|_{r}\nonumber\\
&&\qquad\quad{}+   T^{1/r-1/2}\sum_{n\ge 0}2^{n/2} \max_{0\le
k\le
2^n-1}\big\| |\xi_{2^n+k}-\widehat
\xi_{2^n+k}^{N_{2^n+k}}|^r\big\|^{1/r}_{1}\nonumber \\
&&\qquad= T^{ 1/r-1/2}\|\xi_0-\widehat
\xi_0^{N_0}\|_{r}\nonumber\\
&&\qquad\quad{}+   T^{ 1/r-1/2}\sum_{n\ge 0}2^{
n/2} \max_{0\le
k\le 2^n-1}\big\|
\xi_{2^n+k}-\widehat
\xi_{2^n+k}^{N_{2^n+k}}\big\|_{r}.\nonumber
\end{eqnarray}
\end{cas}

Let $\delta:=\rho-r$. It follows from Lemma~\ref{UPierce} (Pierce
lemma) that, for every
$N\ge 1$ and every r.v.  $\xi\!\in L^r(\mathbb{P})$,

%e2.12 ###
\begin{equation}\label{Piercexi}
\inf_{{\operatorname{card}}(\alpha)\le N}\|\xi -\widehat \xi^{\alpha}\|_{r}
\le C_{r,\rho}
\|\xi\|_{\rho}N^{-1}.
\end{equation}

Now, using the monotony in $p$ of the $L^p$-norms with respect to the
probability measure $2^{n+1}\mbox{\bf
1}_{[0,2^{-(n+1)}T]}(t)\,dt/T$, Fubini's theorem, the $(L_{r,
\varphi})$-Lipschitz continuity assumption~(\ref{assumption1})(i)
and (\ref{1.2}), we obtain
%e2.13 ###
\begin{eqnarray}\label{Fonda0}
\qquad&&\mathbb{E}\,|\xi_{2^n+k}|^{\rho}\nonumber\\
&&\qquad= \mathbb{E}\,|(X|e_{2^n+k})|^{\rho}\nonumber\\
&&\qquad\le 2^{(
n/2)\rho}T^{-\rho/2}\nonumber\\
&&\qquad\quad{}\times\mathbb{E}\biggl(\int_0^{2^{-(n+1)}T}\big|X_{({k}/{2^n})T+u}-X_{({2k+1})/({2^{n+1}})T+u}\big|\,du\biggr)^{\rho}\nonumber\\
&&\qquad\le  2^{(n/2)\rho}2^{-(n+1)\rho}T^{\rho/2}
\vid
&&\qquad\quad{}\times\mathbb{E}\biggl(\int_0^{2^{-(n+1)}T}\big|X_{({k}/{2^n})T+u}-X_{({2k+1})/({2^{n+1}})T+u}\big|^{\rho}\,2^{n+1}\,du/T\biggr)\nonumber\\
&&\qquad\le 2^{-\rho}2^{-(n/2)\rho+n+1}T^{\rho/2-1}\nonumber\\
&&\qquad\quad{}\times\int_0^{2^{-(n+1)}T}\mathbb{E}\big|X_{
({k}/{2^n})T+u}-X_{({2k+1})/({2^{n+1}})T+u}\big|^{\rho}\,du \nonumber \\
&&\qquad\le  2^{-(n/2)\rho+n+1-\rho}T^{\rho/2-1}
\int_0^{2^{-(n+1)}T}\bigl(\varphi(T/2^{n+1})\bigr)^{\rho}\, du\nonumber\\
&&\qquad\le C_{X,T,r,\rho}2^{-(n/2) \rho}\bigl(\varphi(T/2^{n+1})\bigr)^\rho.\nonumber
\end{eqnarray}

At this stage, we    assume a priori that  the size
sequence\break    $(N_{2^n+k})_{n\ge 0,\, k=0,\ldots, 2^{n-1}}$ of the
marginal codebooks is nonincreasing as $2^n+k$  increases and
satisfies
\[
1\le  \prod_{k\ge 0}  N_{k} \le N.
\]

We assume that all the quantizations induced by these  codebooks
are   $L^r$-optimal up to  $n\le m$, that is,
\begin{eqnarray*}
&&\|\xi_{2^n+k}-\widehat \xi_{2^n+k}\|_{r}\\
&&\qquad=\inf_{\operatorname{card}(\alpha)\le
N_{2^n+k}} \|\xi_{2^n+k}-\widehat
\xi_{2^n+k}^\alpha\|_{r}
\end{eqnarray*}
and that  $\widehat
\xi_{2^n+k}=0$ otherwise. Then,
 combining~(\ref{Fonda1}),~(\ref{Fonda0}) and~(\ref{Piercexi}) (Pierce
Lemma) yields
\begin{eqnarray*}
\big\||X-\widehat X|_{L^r_{T}}\big\|_{r}
&\le& C_{X,T,
r,\rho}\Biggl(\frac{1}{N_0} +\sum_{n\ge 0}
\frac{\varphi(T2^{-(n+1)})}{N_{2^{n+1}}}\Biggr)\\
&\le&  C_{X,T,
r,\rho}\Biggl(\frac{1}{N_0} +\frac 1T\sum_{n\ge
0}\sum_{k=0}^{2^{n+1}-1}
\frac{\Phi(2T/(2^{n+1}+k))}{N_{2^{n+1}+k}}\Biggr)\\
&=& C_{X,T, r,\rho}\Biggl(\frac{1}{N_0} +\frac 1T\sum_{k\ge 2}
\frac{\Phi(2T/k)}{N_{k}}\Biggl),
\end{eqnarray*}
where $\Phi(x):= x\varphi(x)$, $x\!\in(0,+\infty)$.  This function
$\Phi$ is regularly varying (at $0$) with index $b+1$. This implies,
in particular, that
  there is a real constant $c>0$ such that $\Phi(T/k)\le
c\Phi(1/(k+1))$ for every $k\ge 2$. Hence, inserting, for
convenience,  the term
$\Phi(1/2)/N_1$ and modifying the real constant $C_{X,T, r,\rho}$ in
an appropriate way finally yields
\[
\big\||X-\widehat X|_{L^r_{T}}\big\|_{r}\le C_{X,T, r,\rho}  \sum_{k\ge
1} \frac{\Phi(1/k)}{N_{k-1}}.
\]
Now, set, for convenience, $\nu_k =\Phi(1/k)$, $k\ge 1$. Note that in
the case $b=0$, the integrability condition $\int_0^1
\varphi(\xi)/\xi\,d\xi<+\infty$ implies $\sum_k \nu_k<+\infty$.
Consequently, an upper bound for the quantization rate is given by the
solution of the following optimal allocation
problem:
%e2.15 ###
%e2.14 ###
\begin{eqnarray}\label{Opti1}
e_{N,r}(X,L^r_{T})&\le&
C_{X,T,r,\rho}
\min\Biggl\{
\sum_{k\ge 1} \frac{\nu_k}{N_{k-1}},\nonumber\\
&&\hspace{61pt}\prod_{k\ge 0}N_k \le N, \;
N_0\ge \cdots\ge N_k \ge \cdots\ge 1
\Biggr\}
\vid
&=& C_{X,T,r,\rho}
\min\Biggl\{ \sum_{k= 1}^m \frac{\nu_{k}}{N_{k-1}}
+\sum_{k\ge m+1}\nu_k,\,m\ge 1,\nonumber\\
&&\hspace{63pt}\prod_{0\le k \le m-1} N_k \le N, \,   N_0\ge
\cdots\ge N_{m-1}\! \ge\!   1
\Biggr\}.\nonumber
\end{eqnarray}

The rest of the proof follows the approach developed in~\cite{LUPA1}
[Section 4.1, especially Lemma 4.2, Theorem~4.6(i)--(iii) and its
proof] and \cite{LUPA2}. However, one must be be aware that we have had to
modify
some notation.\noqed
\end{pf}
\begin{proposition}\label{TechNu}
Assume
$\nu_k= \Phi(1/k)$, $k\ge 1$,  where $\Phi(x)= x\varphi(x)$,\break
$\varphi\dvtx(0,+\infty)$ is
a nondecreasing,  regularly varying function at
$0$ with index $b\ge 0$ with
$\int_{0}^1\varphi(\xi)\frac{d\xi}{\xi}<+\infty$ when $b=0$. Then:
\begin{longlist}
\item $\lim_k \nu_k/\nu_{k+1}=1$;
\item $  (\prod_{k=1}^n
\nu_k)^{ 1/n} \sim e^{b+1} \nu_n$;
\item  $\sum_{k=n+1}^{\infty} \nu_k + n
\nu_k \sim c\psi(1/n)$, where $c =1+1/b$ if $b>0$; $c=1$ if
$b=0$;
\end{longlist}
\[
\psi(x) = \varphi(x) \qquad \mbox{ if } b>0;\qquad \psi(x) :=
\int_0^x\varphi(\xi)\frac{d\xi}{\xi}
\qquad \mbox{ if }b=0.
\]
(See \textup{\cite{LUPA1}} for a proof.)
\end{proposition}
\begin{pf*}{Proof of Theorem \protect\ref{Upper1} (C\normalfont{ontinued})}
Set
\[
m=m^*(N)= \max\Biggl\{m\ge 1 \dvtx N^{
1/m} \nu_m \Biggl(\prod_{j=1}^m
\nu_j\Biggr)^{-1/m}\ge 1\Biggr\}
\]
and
\[
N_{k-1} = N_{k-1}(N):= \Biggl[N^{1/m} \nu_k
\Biggl(\prod_{j=1}^m
\nu_j\Biggr)^{-1/m}\Biggr]\ge 1,\qquad k=1,\ldots,m.
\]
It follows from  Proposition~\ref{TechNu}(ii) that
\[
m=m^*(N)\sim \frac{\log N}{b+1} \qquad \mbox{as }  N\to \infty.
\]
Then
\begin{eqnarray*}
\sum_{k= 1}^m
\frac{\nu_{k}}{N_{k-1}}&\le &\max_{k\ge 1}(1+1/N_{k-1})mN^{-
1/m}
\Biggl(\prod_{j=1}^m
\nu_j\Biggr)^{1/m} \\
&\le& 2mN^{-1/m}
  \Biggl(\prod_{j=1}^m \nu_j\Biggr)^{1/m}\\
&\le&  2 m\nu_m.
\end{eqnarray*}
Consequently, this time  using (iii) in Proposition~\ref{TechNu},
\begin{eqnarray*}
\sum_{k= 1}^m \frac{\nu_{k}}{N_{k-1}}+\sum_{k\ge
m+1}\nu_k&\le& 2\Biggl(m\nu_m +\sum_{k\ge m+1 } \nu_k\Biggr)\\
&=& O\bigl(\psi(1/\log N)\bigr)
\end{eqnarray*}
so that
\[
  \big\||X-\widehat X|_{L^p_{T}}\big\|_{r}= O\bigl(\psi(1/\log N)\bigr).
\]

\begin{cas}[($\rho\le 1$)]
Here, we rely on the
pseudo-triangular inequality
\[
|f+g|^r_{L^r_{T}} \le |f|^r_{L^r_{T}}+|g|^r_{L^r_{T}},
\]
which follows from the elementary inequality $(u+v)^r\le u^r+v^r$:
\begin{eqnarray*}
|X-\widehat X|^r_{L^r_{T}}&\le& T^{1- r/2}|\xi_0-\widehat
\xi_0^{N_0}|^r + \sum_{n\ge 0}
\Bigg|\sum_{k=0}^{2^n-1}(\xi_{2^n+k}-\widehat
\xi_{2^n+k}^{N_{2^n+k}})e_{2^n+k}\Bigg|^r_{L^r_{T}}\\
&=& T^{1-r/2}|\xi_0-\widehat \xi_0^{N_0}|^r + T^{1-
r/2}\sum_{n\ge 0}2^{n({r}/{2} - 1)}
\sum_{k=0}^{2^n-1}|\xi_{2^n+k}-\widehat
\xi_{2^n+k}^{N_{2^n+k}}|^r
\end{eqnarray*}
so that
%e2.16 ###
\begin{eqnarray}\label{Fonda1prime}
\qquad\quad\big\||X-\widehat X|_{L^r_{T}}
\big\|^r_{r}&=&\big\||X-\widehat X|^r_{L^r_{T}}\big\|_{1}\nonumber\\
&\le& T^{1- r/2}\big\||\xi_0-\widehat \xi_0^{N_0}|^r \big\|_{1}\nonumber\\
&&{}+T^{1- r/2}\sum_{n\ge 0}2^{n(r/2
-1)}\Bigg\|
\sum_{k=0}^{2^n-1}|\xi_{2^n+k}-\widehat
\xi_{2^n+k}^{N_{2^n+k}}|^r\Bigg\|_{1}\nonumber\\
&\le & T^{1- r/2}\|\xi_0-\widehat \xi_0^{N_0}\|^r_{r}\nonumber\\
&&{}+T^{1-r/2}\sum_{n\ge 0}2^{n( r/2 -1)}  2^n
\!\!\!\max_{0\le
k\le 2^n-1}\big\|
|\xi_{2^n+k}-\widehat
\xi_{2^n+k}^{N_{2^n+k}}|^r\big\|_{1}\\
&=& T^{1-r/2}\|\xi_0-\widehat \xi_0^{N_0}\|^r_{r}\nonumber\\
&&{}+T^{1-r/2}\sum_{n\ge 0}2^{{nr}/{2}} \max_{0\le
k\le
2^n-1}\big\||\xi_{2^n+k}-\widehat
\xi_{2^n+k}^{N_{2^n+k}}|^r\big\|_{1}\nonumber\\
&=& T^{1-r/2} \|\xi_0-\widehat
\xi_0^{N_0}\|^r_{r}\nonumber\\
&&{}+ T^{1-r/2}\sum_{n\ge 0}2^{{nr}/{2}}
\max_{0\le
k\le
2^n-1}\|
\xi_{2^n+k}-\widehat
\xi_{2^n+k}^{N_{2^n+k}}\|^r_{r}.\nonumber
\end{eqnarray}
This inequality replaces~(\ref{Fonda1}). We then note that
\begin{eqnarray*}
\mathbb{E}|\xi_{2^n+k}|^{\rho}
&\le& 2^{
(n/2)\rho}T^{-\rho/2}\bigl(2^{-(n+1)}T
\varphi(T/2^{n+1}) \bigr)^{\rho}\\
&= &
C_{X,T,r,\rho}2^{-(n/2) \rho}\bigl( \varphi(T/2^{n+1})\bigr)^\rho
\end{eqnarray*}
so that
\[
\big\||X-\widehat
X|_{L^r_{T}}\big\|_{r}^r\le  C_{X,T,r,\rho}\Biggl(\frac{1}{N_0^r}+
\sum_{n\ge 0}
\frac{\varphi(T2^{-(n+1)})^r}{N^r_{2^{n+1}}}\Biggr).
\]
We then set $\widetilde \varphi(u) = (\varphi(u))^r$, $\widetilde
N_k = N_k^r$ and $\widetilde N:= N^r$. We proceed
for
$\||X-\widehat X|_{L^r_{T}}\|_{r}^r$ with these ``tilded''
parameters as for $\||X-\widehat
X|_{L^r_{T}}\|_{r}$    in
the case
$\rho>1$.\quad\qed
\end{cas}\noqed
\end{pf*}
\begin{remark*}
\textit{Concerning the case} $p>r$.  When  $p\ge \rho > r $, the
$(L^r(\mathbb{P}),L^p_{T})$-quantization problem remains consistent.
However,
there is  a price to be paid for considering  a $p$  exponent
greater than $\rho$. Thus, if $\varphi$ in $(L_{(\rho,\varphi)})$
has
regular variations with exponents $b>0$ at $0$ and if
$b+\frac
1p-\frac 1r>0$, then the same approach yields the rate
\[
e_{N,r}(X,L^p_{T})\le C_{X,r,\delta,T,p} \varphi(1/\log N)(\log
N)^{1/r-1/p}.
\]
We do not know whether  it is due to
our approach or if it is the best possible universal  rate.

\textit{Concerning lower bounds}. In several
situations, when the   assumption $(L_{\rho,\varphi})$ is
optimal in
terms of mean regularity of a process, the upper bound for the
functional quantization rate turns out to be the true rate.
We have
no general
result in that direction so far since most lower bound
results rely on  a different approach, namely the small
deviation theory.
Thus,
 in \cite{GRLUPA1}, a connection is
established between (functional) quantization and small
deviation for Gaussian processes.
In particular, this approach
provides a method to derive a lower bound for the $(L^r(\mathbb{P}),
L^p_{T})$-quantization rate from some upper bound for the small
ball
problem. A careful reading of the paper (see the proof of Theorem~1.2
in~\cite{GRLUPA1}) shows that this small deviation lower bound holds
for any
\textit{unimodal} (w.r.t. 0) nonzero process. To be precise,
let $p\!\in (0,\infty)$ and assume that  $\mathbb{P}_{X}$ is
$L^p_{T}$-unimodal in the following
sense: there exists a real
$\varepsilon_0 > 0$ such that
\[
\forall x\!\in
L^p_{T},\forall\, \varepsilon \in (0,\varepsilon_0]\qquad
\mathbb{P}(|X-x|_{L^p_{T}}\le \varepsilon)\le \mathbb{P}(|X|_{L^p_{T}}\le
\varepsilon).
\]
(For centered Gaussian processes, this follows for
$p\!\ge\! 1$ from  Anderson's inequality.) If

\[
G\bigl(-\log\bigl(\mathbb{P}(|X|_{L^p_{T}}\le \varepsilon)\bigr)\bigr)=
\Omega(1/\varepsilon) \qquad \mbox{as }\varepsilon\to 0
\]
for
some increasing unbounded function $G\dvtx(0,\infty)\to (0,\infty)$,
then
%e2.17 ###
\begin{equation}\label{sblowerbound}
\forall r\!\in
(0,\infty), \forall c>1\qquad e_{{N,r}}(X, L^p_{T}) =
\Omega\biggl(\frac{1}{G(\log(cN))}\biggr).
\end{equation}
\end{remark*}

%s3 ###
\section{Applications and examples}\label{FirstAppli}
In this
section, we give some examples which illustrate that the upper bound
derived from the mean pathwise regularity may be
optimal or
not.

%s3.1 ###
\subsection{Application to It\^o processes and  $d$-dimensional
diffusion processes}
\mbox{}

Let $W$ denote  an $\mathbb{R}^d$-valued
standard Brownian motion defined on a probability space
$(\Omega,\mathcal{A},\mathbb{P})$ and let
$(\mathcal{F}_t^W)_{t\in [0,T]}$ denote
its natural filtration (completed with all the $\mathbb{P}$-negligible sets).
Let $X$ be a
$1$-dimensional It\^o process   defined by
\[
dX_t =
G_t\,dt+H_t \cdot dW_t, \qquad X_0= x_0\!\in \mathbb{R},
\]
where $(G_t)_{t\in
[0,T]}$ is a  real-valued  process and $(H_t)_{t \in [0,T]}$ is an
$\mathbb{R}^d$-valued   process, both  assumed $(\mathcal{F}_t^W)_{t\in
[0,T]}$-progressively measurable.  Assume that there exists a real
number
$\rho\ge 2$ such that

%e3.1 ###
\begin{equation}\label{ItoAssump}
\sup_{t\in[0,T]}\mathbb{E}|G_t|^\rho
+\sup_{t\in[0,T]}\mathbb{E}|H_t|^\rho<+\infty,
\end{equation}
where $|\cdot|$
denotes any norm on $\mathbb{R}^d$. Then
(see,
e.g.,~\cite{BOLE}) the $\varphi$-Lipschitz assumption
$(L_{\varphi,\rho})$(i)  [i.e., (\ref{assumption1})(i)] is
satisfied  with $\varphi(u)= cu^{1/2}$.
It follows
from~Theorem~\ref{Upper1} that
\[
\forall\, r,\,p\!\in (0,\rho)\qquad
e_{N,r}(X,L^p_{T})  = O((\log
N)^{-1/2}).
\]

Let
$X=(X^1,\ldots,X^d)$ be an $\mathbb{R}^d$-valued diffusion process  defined by
\[
dX_t = b(t,X_t)\,dt+\sigma(t,X_t) \, dW_t, \qquad X_0= x_0\!\in \mathbb{R}^d,
\]
where $b\dvtx[0,T]\times \mathbb{R}^d\to \mathbb{R}^d$ and $\sigma\dvtx[0,T]
\times\mathbb{R}^d\to \mathcal{M}(d\times q,\mathbb{R})$ are  Borel functions satisfying
\[
\forall\, t\in[0,T],\  \forall x\in\mathbb{R}^d \qquad
|b(t,x)|+\|\sigma(t,x)\| \le C(1+|x|)
\]
and
$W$ is an
$\mathbb{R}^q$-valued  standard Brownian motion. The above
assumption does not imply that such a diffusion process $X$ exists.
(The
existence holds provided $b$ and $\sigma$ are Lipschitz in $x$
uniformly with respect to $t\!\in[0,T]$.) Then, every component
$X^i$ is an It\^o process [with
$G_t=b^i(t,X_t)$ and
$H_t:=
\sigma^{i\bolds{\cdot}}(t,X_t)$] for which assumption~(\ref{ItoAssump})  is
satisfied for every $\rho>0$ (see, e.g.,~\cite{BOLE}). On the other
hand, if
$(u^1,\ldots,u^d)$ denotes the canonical basis of $\mathbb{R}^d$ and
$|\cdot|$ denotes any norm on $\mathbb{R}^d$, then for every $p\ge 1$
and every
$f:=\sum_{1\le i\le d} f^iu^i\dvtx[0,T]\to
\mathbb{R}^d$,
\[
|f|_{L^p_{\mathbb{R}^d}([0,T],dt)}\le \sum_{i=1}^d |f^i|_{L^p_T}|u^i|.
\]
Now, we can quantize each It\^o process $(X_t^i)_{t\in[0,T]}$,
$i=1,\ldots,d,$ using an $(L^r,L^p_{T})$-optimal quantizer
$\alpha^{(i)}$ of size $[\sqrt[d]{N}]$. It is clear that the
resulting product quantizer $\prod_{i=1}^d \alpha^{(i)}$ of
size $[\sqrt[d]{N}]^d\le N$ induces an
$(L^r,L^p_{\mathbb{R}^d}([0,T],dt))$-quantization error $O((\log
N)^{-1/2})$ (see,
e.g.,~\cite{LUPA4}). Combining these obvious remarks finally yields
\[
\forall\, r,\,p>0 \qquad e_{N,r}(X,L^p_{\mathbb{R}^d}([0,T],dt))=
O((\log N)^{-1/2}).
\]
In the ``smooth'' case $H\equiv 0$, the regularity assumption
$(L_{\varphi,\rho})$ is satisfied with $\varphi(u)= cu$. We obtain the
universal upper bound
\[
\forall\, r,\,p\!\in (0,\rho)\qquad e_{{N,r}}(X,L^p_{T}) = O((\log N)^{-1}).
\]
Both rates are optimal as universal rates for $p\ge 1$, as can be seen
from $X=W$ and $X= \int_0^{\bolds{\cdot}}G_s\,ds$ with
$G_t= \int_0^t (t-s)^{\beta-1/2}\,dW_s$ ($\beta>0$
and $d=1$), respectively (see~\cite{GRLUPA1}).

As far as quantization rates are concerned, this extends to general
$d$-dimensional diffusions a result obtained
in~\cite{LUPA4} by stochastic calculus techniques for a more
restricted class of Brownian diffusions (which
includes $1$-dimensional ones). This also extends (the upper bound
part of the) the result obtained
in~\cite{Dereich} for another class of (essentially $1$-dimensional)
Brownian diffusions. For the class
investigated in~\cite{LUPA4}, it is shown that   under an ellipticity
assumption  on
$\sigma$,  this rate is optimal in the case $r,p\ge 1$.
In~\cite{Dereich}, still with a (mild)  ellipticity assumption, the
rate  is sharp for $p\ge 1$.   This leads us to
conjecture that this rate is optimal for not too degenerate Brownian
diffusions.

%s3.2 ###
\subsection{Application to fractional Brownian motion}
The fractional Brownian motion $W^H$ with Hurst constant
$H\!\in(0,1]$  is a Gaussian process satisfying, for every $\rho>0$,
\[
\mathbb{E}|W^H_t-W^H_s|^\rho = C_{H,p} |t-s|^{\rho H}\quad \mbox{and}
\quad (W^H_s)_{0\le s\le t} \stackrel{\mathcal{L}}{\sim}t^H
(W^H_{s/t})_{0\le s\le
t}.
\]
So, using Theorem~\ref{Upper1}, we obtain  $e_{N,r}(W^H,
L^p_{T})=O((\log N)^{-H})$ as an
$(L^r(\mathbb{P}),\break|\cdot|_{L^p_{T}})$-quantization  rate for every $r,p>0$.
This rate is known to be optimal for $p\ge 1$. In fact, a sharp  rate is
established  (see~\cite{LUPA2}, when
$p=r=2$, or~\cite{Dereichetal}) [i.e., the computation of the exact
value of $\lim_N N(\log N)^H e_{{N,r}}(W^H,L^p_{T})$].

%s3.3 ###
\subsection{Stationary processes}
Let $X$ be a centered weakly (square-integrable) stationary process.
Then
\[
\mathbb{E }|X_t-X_s|^2 = \mathbb{E } |X_{t-s}-X_0|^2 =2\,{\operatorname{Var}}(X_0)\bigl(1-
c(|t-s|)\bigr),
\]
where $c(t)$ denotes the correlation between $X_t$ and $X_0$. Hence, if
\[
c(u)= 1- \kappa u^{2a}+o(u^{2a})\qquad \mbox{as }u\to 0,
\]
then the $L^r(\mathbb{P})$-rate for $L^p_{T}$-quantization $0<p,r< 2$, will
be  $O((\log(N))^{-a})$.  If, furthermore, $X$ is  a Gaussian
process (like the Ornstein--Uhlenbeck process with $a=1/2$), then
this $O((\log N)^{-a})$ rate holds for any $r,\,p>0$ since, for every
$\rho\!\in \mathbb{N}^*$,
\[
\mathbb{E} |X_t-X_s|^\rho = \mathbb{E} |X_{t-s}-X_0|^\rho =C_\rho \bigl({\operatorname{ Var}}(X_0)\bigl(1-
c(|t-s|)\bigr)\bigr)^{\rho/2}.
\]

%s3.4 ###
\subsection{Self-similar processes with stationary increments}
Let $X=(X_t)_{t\in [0,T]}$ be an $H$-self-similar process with
stationary increments [$H\!\in(0,\infty)$]. Assume $X_1\!\in
L^\rho(\mathbb{P})$ for some $\rho\ge 1$. Then
\[
\mathbb{E} |X_t-X_s|^\rho = C_\rho |t-s|^{\rho H}
\]
for every $s,\, t\!\in [0,T]$. Since $X$ is stochastically
continuous, it has a   bimeasurable modification.
Theorem~\ref{Upper1} then gives
\[
\forall\, r,\, p \!\in(0,\rho) \qquad e_{N,r}(X,L^p_{T}) = O( (\log N)^{-H}).
\]

If, furthermore, $X$ is $\alpha$-stable, $\alpha\!\in(1,2)$, then
$X_1 \!\in L^\rho(\mathbb{P})$ for every $\rho \!\in [1,\alpha)$ so that
\[
\forall\, r,\, p \!\in(0,\alpha) \qquad e_{N,r}(X,L^p_{T}) = O(
(\log N)^{-H}).
\]
This class of examples comprises, for example, the linear $H$-fractional
$\alpha$-motions with $\alpha \!\in(1,2)$, $H\!\in(0,1)$ and the
$\log$-fractional $\alpha$-stable motions
with $\alpha \!\in(1,2)$, where $H= 1/\alpha$ (see~\cite{EMMA,SATA}).

%s3.5 ###
\subsection{L\'evy processes\textup{:} a first approach}\label{Levyproc1}
A (c\`adl\`ag) L\'evy process $X=\break (X_t)_{t\in \mathbb{R}_+}$---or Process with
Stationary
Independent Increments ($\mathit{PSII}$)---is characterized by its so-called
\textit{local characteristics} appearing in the L\'evy--Khintchine formula (for
an
introduction to L\'evy processes, we refer to~\cite{BER,JASH,SAT}).
These characteristics depend on the way the ``big'' jumps are
truncated. We will adopt, in the following, the
convention that the truncation occurs at   size
1. So that, for every $t\!\in \mathbb{R}_+$,
%%%%%%%%%%%%%%%%%%%%%%%%%%%%%%%%%%%%%%%%%%%%%%%%%%%%%%%%%%%%%%%%%%%%%%%rt
\begin{eqnarray*}
&&\mathbb{E}(e^{iuX_t})= e^{-t\psi(u)}\\
&&\qquad\mbox{where }\psi(u)=
-iua+\tfrac {1}{2} \sigma^2u^2 -
\int_{\mathbb{R}\setminus\{0\}}\bigl(e^{iux}-1-iux\mathbf{1}_{\{|x|\le
1\}}\bigr)\nu(dx),
\end{eqnarray*}
where $a, \sigma\!\in \mathbb{R}$ and  $\nu$ is a nonnegative measure on
$\mathbb{R}\setminus\{0\}$ such that $\nu(x^2\wedge 1)<+\infty$. The measure
$\nu$ is called
the \textit{L\'evy measure} of the process. It can be shown that a L\'evy
process is a compound Poisson process if and only if $\nu$ is a
finite measure and  has finite variation if and only if
$
\int_{\{|x|\le 1\}} |x|\nu(dx)<+\infty$. Furthermore,
\[
X_t \!\in L^\rho(\mathbb{P})\quad \mbox{if and only if }\quad
\int_{\{|x|\ge 1\}} |x|^\rho\nu(dx)<+\infty.
\]

We will  extensively use the following Compensation Formula (see,
e.g.,~\cite{BER} page 7):
%e3.2 ###
\begin{equation}\label{compensation}
\qquad\quad\mathbb{E} \sum_{s\ge 0} F(s,X_{s-},\Delta X_s)\mathbf{1}_{\{\Delta X_s\neq
0\}} =
\mathbb{E}
\int_{\mathbb{R}_+}ds\int_{\mathbb{R}\setminus\{0\}} F(s,
X_{s-},\xi)\nu(d\xi),
\end{equation}
where $F\dvtx \mathbb{R}_+\times \mathbb{R}^2\to \overline \mathbb{R}_+$ is  a Borel function. As
concerns assumption~(\ref{assumption1}), note that  the very
definition of a
L\'evy process implies that
\[
\mathbb{E}|X_t-X_s|^\rho = \mathbb{E}|X_{t-s}|^\rho \quad\mbox{and}\quad
\mathbb{E}\sup_{ s\in[t,t+h]}|X_t-X_s|^\rho = \mathbb{E}\sup_{  s\in[0, h]}
|X_{s}|^\rho,
\]
so we may focus on the distribution of  $X_t$ and $X^*_t:
=\sup_{s\in[0,t]}|X_s|$. Finally, note that it follows from the usual
symmetry principle
(see~\cite{SAT}) that for any L\'evy process, $\mathbb{P}(X^*_t>u+v) \le
\mathbb{P}(|X_t|>u)/\mathbb{P}(X^*_t
\le v/2)$ so that $\mathbb{E}|X_t|^r$ and $\mathbb{E}|X^*_t|^r$ are simultaneously
finite or infinite when $r>0$.

The following result  is established  in~\cite{MI}.
\begin{lema}[(Millar's Lemma)]\label{Millar0}
 Assume $\sigma=0$. If
there exists a real number
$\rho\!\in(0,2]$ such that $
\int_{\mathbb{R}\setminus\{0\}}|x|^\rho
\nu(dx)<+\infty $, then there exist some real constants
$a_\rho\!\in\mathbb{R}$ and  $C_\rho>0$ such that
%e3.3 ###
\begin{equation}\label{Millar}
\forall\, t\!\ge 0\qquad   \mathbb{E}\biggl(\sup_{s\in[0,t]}|X_s-a_\rho
s|^\rho\biggr)\le C_\rho t.
\end{equation}
Furthermore, one may set $a_\rho=0$ if $\rho\ge 1$.
\end{lema}

Hence, it follows as a  consequence of Theorem~\ref{Upper1} that
%e3.4 ###
\begin{equation}\label{MillarRate}
\forall\, r,\,p\!\in (0,\rho)\qquad e_{N,r}(X,L^p_{T}) =O((\log
N)^{-1/\rho}).
\end{equation}

This follows from the  following straightforward  remark:
if  $\beta\subset L^p_{T}$ is   an\break  $N$-quantizer and
$\xi\!\in L^p_{T}$ [here $\xi(t)= a_\rho t$], then
\[
\big\||X-\widehat
X^\beta|_{L^p_{T}}\big\|_{r}=\big\||(X+\xi)-\widehat{(X+\xi)}^{\xi+\beta}|_{L^p_{T}}\big\|_{r}\qquad
\mbox{with }  \xi+\beta= \{ \xi+f,\;f\!\in\beta\}.
\]

However, rate~(\ref{MillarRate})  may be suboptimal,  as illustrated
below with $\alpha$-stable processes and Poisson processes.
In Section~\ref{LevyLevy}, we establish two improvements of this rate
under some natural hypotheses (see~Theorem~\ref{Levyrate}
for a broad class of   L\'evy processes with infinite L\'evy
measure  and~Proposition~\ref{Poisson} for compound Poisson
processes).

\textit{The
$\alpha$-stable processes}.\label{stable1}
The (strictly) $\alpha$-stable processes are families of L\'evy
processes indexed by $\alpha\!\in(0,2) $ satisfying  a
self-similarity property, namely
\[
\forall\, t\!\in\mathbb{R}_+\qquad X_t \mathop{\sim}^\mathcal{L}  t^{
1/\alpha}
X_1\quad\mbox{and}\quad\sup_{0\le s\le t} |X_s|
\mathop{\sim}^\mathcal{L}  t^{ 1/\alpha} \sup_{0\le s\le 1}|X_s|.
\]
Furthermore,
\[
\sup\biggl\{r \dvtx \mathbb{E}\biggl(\sup_{0\le s\le 1}|X_s|^r\biggr) <+\infty\biggr\} =
\alpha\quad
\mbox{and}\quad \mathbb{E }|X_1|^\alpha =+\infty.
\]

Consequently, it follows from Theorem~\ref{Upper1}, applied with
$\varphi(u) := u^{{1}/{\alpha}}$, that
%
%e3.5 ###
\begin{equation}\label{stablerate}
\forall\, p,\,r\!\in(0,\alpha)\qquad e_{N,r}(X, L^p_{T})
= O\biggl(\frac{1}{(\log N)^{{1}/{\alpha}}}\biggr).
\end{equation}

In the symmetric case, an $\alpha$-stable process $X$ being
subordinated to a Brownian motion ($X_t= W_{A_t}$ with $A$ a one-sided
$\alpha/2$-stable process) has a unimodal distribution by the
Anderson inequality (see Section~\ref{suboLevy} below, entirely
devoted to subordinated L\'evy processes).  Substituting
into~(\ref{sblowerbound})   the small deviation estimates established
in  \cite{LISH}    shows the rate optimality   of our upper bound for
$e_{{N,r}}$ when $p\ge 1$, that is,
%e3.6 ###
\begin{equation}\label{alphasblb}
\forall\, r\!\in (0,\alpha),\, \forall\, p\!\in [1,\alpha)\qquad
e_{{N,r}}(X,L^p_{T})\approx (\log N)^{- 1/\alpha}.
\end{equation}

\textit{The $\Gamma$-processes}.
These are subordinators (nondecreasing L\'evy processes) whose
distribution $\mathbb{P}_{X_t}$ at time $t$ is a
$\gamma(\alpha,t)$-distribution,
\[
\mathbb{P}_{X_t}(dx)= \frac{\alpha^t}{\Gamma(t)}\mathbf{1}_{(0,\infty)}(x)
x^{t-1}e^{-\alpha x} \,dx.
\]
So, easy computations show that for every $\rho>0$,
\[
\mathbb{E}|X_t|^{\rho}=  \frac{\Gamma(t+\rho)}{\alpha^\rho
\Gamma(t+1)}\,t\sim \frac{\Gamma(\rho)}{\alpha^\rho
\Gamma(1)}\,t\qquad \mbox{as  }
t\to 0.
\]
Consequently, it follows from Theorem~\ref{Upper1} that
\[
\forall\, p\!\in(0,+\infty), \, \forall\,r\,\in(0, p] \qquad
e_{N,r}(X, L^p_{T})=
O\biggl(\frac{1}{(\log (N))^{1/p-\varepsilon}}\biggr)\qquad
\forall\, \varepsilon>0.
\]

\textit{Compound Poisson processes from the mean
regularity viewpoint}. One considers
a compound Poisson process
\[
  X_t= \sum_{k=1}^{K_t} U_k,
\]
where $K=(K_t)_{t\in[0,T]}$ denotes a standard Poisson process with
intensity $\lambda=1$ defined on a probability space
$(\Omega,\mathcal{A}, \mathbb{P})$ and
$(U_k)_{k\ge 1}$ an  i.i.d. sequence of random variables defined on
the same probability space, with $U_1\!\in L^{\rho}(\mathbb{P})$ for some
$\rho>0$.  Then,
standard computations show that
%e3.8 ###
%e3.7 ###
\begin{eqnarray}
\label{Poissrle1}
\qquad\mathbb{E}\sup_{0\le s\le t}\Bigg| \sum_{k=1}^{K_s} U_k
\Bigg|^\rho&\le& \mathbb{E}\sum_{k=1}^{K_t} |U_k|^\rho =
t\,\|U_1\|_{\rho}^\rho \qquad\hspace{28pt}\mbox{if } 0<\rho\le 1,\\
\label{Poissrge1}\mathbb{E}\Bigg| \sum_{k=1}^{K_t} U_k \Bigg|^\rho&\le&
t\|U_1\|_{\rho}^\rho \times \Biggl[e^{-t}\sum_{k\ge
1}\frac{t^{k-1}k^\rho}{k!}\Biggr]\qquad
\mbox{if } \rho> 1.
\end{eqnarray}
Consequently, assumption~(\ref{assumption1}) is
fulfilled with $\varphi(u)=cu^b$, where
$b=1/\rho$ and $c$ is a positive real constant. Theorem~\ref{upper-bound} then yields
\[
\forall\, r,\,p\!\in (0,\rho) \qquad e_{N,r}(X,L^p_{T}) = O((\log
N)^{-1/ \rho}).
\]
Note that when $\rho \le 2$, this is a special case
of~(\ref{Millar}). These rates are  very far from  optimality,
as will be seen further on (in Section~\ref{LevyLevy},  some
faster rates are established by a completely different approach
based on the almost finite-dimensional feature of the paths of such
elementary jump processes). This will emphasize the fact that the mean
regularity of
$t\mapsto X_t$ does not always  control  the quantization rate.

%s4 ###
\section{A  quantization rate for general  L\'evy processes without
Brownian component}\label{LevyLevy}
The aim of this section is to provide a general  result for L\'evy
processes without Brownian component,  with special attention being paid
to  compound Poisson processes
which  appear as a critical case of the main theorem. Before stating
the main results,  we need some further notation related
to L\'evy processes. Set
%e4.1 ###
\begin{eqnarray}
\underline \theta &:=& \inf\biggl\{\theta >0\dvtx \int_{\{|x|\le
1\}}|x|^\theta \nu(dx)<+\infty\biggr\}\!\in [0,2],\label{thetabas}\\
r^* &:=& \sup\biggl\{r>0 \dvtx \int_{\{|x|> 1\}}|x|^r
\nu(dx)<+\infty\biggr\}\le +\infty.\label{thetahaut}
\end{eqnarray}

The exponent $\underline \theta$ is known as the \textit{Blumenthal--Getoor
index} of   $X$ [and is often denoted $\beta(X)$ in the literature].
We define on
$(0,\infty)$ the tail function of the L\'evy measure
$\nu\dvtx u\mapsto \underline
\nu(u):=\nu([-u,u]^c)$. Finally, we set, for every  $\underline
\theta>0$,  $\underline \ell(t) :=  t\underline
\nu(t^{{1}/{\underline\theta}})$ and, for every $\rho>0$,
\[
\Lambda_\rho(t) :=
(\underline\ell(t))^{{1}/{2}}+ (\underline \ell(t))^{
{1}/{\rho}}+ (\underline \ell(t))^{
{2}/{\rho}}\mathbf{1}_{\underline \theta\in (1,2]\cup IV(1)},
\]
where $IV(1)= \varnothing$ if $\underline \theta=1$ and
$\nu(|x|)<+\infty$, and  $IV(1) = \{1\}$ if
$\underline \theta=1$ and $\nu(|x|)=+\infty$.
\begin{teorema}\label{Levyrate}
Let $X=(X_t)_{t\in [0,T]}$ be a L\'evy
process with L\'evy measure
$\nu$ and without Brownian component. Assume
$r^*,\,
\underline
\theta >0$.

\begin{longlist}[(a)]
\item[(a)]  Assume $\underline \theta
\!\in(0,2]\setminus\{1\}$. If $ \int_{\{|x|\le 1\}}
|x|^{\underline \theta} \nu(dx)
<+\infty$ \textup{(}i.e., $\underline
\theta$ holds as a minimum\textup{)} or if the L\'evy measure satisfies
%e4.3 ###
\begin{equation}\label{HypoLevy}
\exists\, c\!\in(0,1], \exists\, C>0\qquad\mathbf{1}_{\{0<|x|\le
c\}}\nu(dx)\le\frac{ C }{|x|^{\underline \theta +1}}
\mathbf{1}_{\{0<|x|\le c\}}\,dx,
\end{equation}
then
%e4.4 ###
\begin{equation}\label{BofRate0}
\forall\, r,\,p\!\in(0,\underline \theta\wedge r^*)\qquad
e_{N,r}(X,L^p_{T})   = O((\log N)^{-{1}/{\underline
\theta}}).
\end{equation}

\item [(b)]Assume $\underline \theta \!\in(0,2)\setminus\{1\}$.
If   the tail function of the L\'evy measure $\nu$ has regular
variation with index
$-b$ at
$0$, then $b=\underline
\theta$ and  the function $\underline \ell$ is
slowly varying   at $0$. If, furthermore, the functions $t\mapsto
t^{{1}/{\underline \theta}}\Lambda_\rho(t)$ are  nondecreasing
in a neighborhood of 0,
then
%
%e4.5 ###
\begin{eqnarray}\label{RateLevy}
\qquad&&\forall\, r,\,p\!\in(0,\underline \theta\wedge r^*) \qquad
e_{N,r}(X,L^p_{T})  =
O((\log N)^{-{1}/{\underline
\theta}}\Lambda_{\rho}((\log\!N)^{-1}))
\vid
&&\qquad\forall\,\rho\!\in(r\vee p, \underline \theta).\nonumber
\end{eqnarray}

\item [(c)]Assume $\underline \theta<r^*$. For every
$r\!\in[\underline \theta,r^*)$ and every $p\in(0,r]$,
%e4.6 ###
\begin{equation}\label{BofRate0b}
    e_{N,r}(X,L^p_{T})   = O((\log N)^{-{1}/{r}+\eta})\qquad\forall\,\eta>0.
  \end{equation}

\item [(d)] When $\underline \theta =1$, if  $\nu$ is symmetric
or $\nu(|x|)<+\infty$, then the above rates~(\ref{BofRate0})
and~(\ref{RateLevy}) are still valid.
\end{longlist}
\end{teorema}

\begin{remark*}
The conclusion in  (a) remains
valid  for any  $\underline \theta\!\in(0,2]$ satisfying
$\int_{\{|x|\le 1\}} |x|^{\underline \theta} \nu(dx)
<+\infty$ or~(\ref{HypoLevy}), not only for the Blumenthal--Getoor
index. In particular, with $\underline \theta\!=\!2$ we obtain
\[
\forall\, r,\,p\!\in(0, 2\wedge r^*)\qquad
e_{N,r}(X,L^p_{T})   = O((\log N)^{-{1}/{2}}).
\]

When $\underline \theta\!\in\{1,2\}$, some rates
can also be derived [even when $\nu$ is not symmetric and
$\nu(|x|)=+\infty$]. Thus, in item~(a), if
$\underline \theta =1$, we can show, by adapting the proof of case
$\underline \theta \!\in(1,2)$ in Proposition~\ref{sous-theta} below,
that
\[
 e_{N,r}(X,L^p_{T})   = O\biggl(\frac{\log\log
N}{\log
N}\biggr).
\]

In most natural settings,
there is  a dominating term in the definition of the function
$\Lambda_\rho$. Thus,
in~(\ref{RateLevy}), we may set
\[
\Lambda_\rho(t)=
\cases{
(\underline \ell(t))^{
{1}/{\rho}}\mathbf{1}_{\{\underline \theta\in (0,1]\setminus
IV(1)\}}+ (\underline \ell(t))^{
{2}/{\rho}}
\mathbf{1}_{{\{\underline \theta\in (1,2]\cup IV(1) \}}},\cr
\hspace{65pt}
\mbox{when $\displaystyle\lim_{t\to 0} \underline \ell(t)=
+\infty$,}\cr
(\underline\ell(t))^{ {1}/{2}},
\qquad \mbox{when  $\displaystyle\lim_{t\to 0} \underline \ell(t)=0$.}
}
\]

Note that this theorem
provides no rate when $\underline \theta =0$, which is the case of an
important class of L\'evy processes
including  compound Poisson
processes. In fact, for these processes, the quantization rate is not
ruled by the mean  regularity of their  paths, as
emphasized in
Section~\ref{CompoundP}.

The proof of this
theorem relies on Theorem~\ref{Upper1}, that is, on the mean pathwise
regularity  of $X$, hence the critical value $\underline \theta$ for
$\rho$
cannot be overcome by such an approach since assumption
$(L_{\varphi,\rho})$ for $\rho>\underline \theta$ would imply that
$X$ has a pathwise continuous modification by the
Kolmogorov criterion.
\end{remark*}
\begin{examples*}
Note that
for \textit{$\alpha$-stable
processes}, $r^*=\underline \theta =
\alpha$, $\nu$ satisfies~(\ref{HypoLevy})  and $\lim_{u\to
0}\underline
\ell(u)
\!\in(0,\infty)$ so that both    rates obtained
from~(\ref{BofRate0}) and~(\ref{RateLevy})   coincide  with
that
obtained in Section~\ref{stable1}, that is, $O((\log
N)^{-{1}/{\alpha}})$. This rate  is most likely optimal.

Let    $\nu^1_{a,\underline \theta}\,(dx): =
\kappa
|x|^{-\underline
\theta -1} (-\log |x|)^{-a}\mbox{\bf
1}_{(0,c]}(|x|)\,dx$, with $0<c<1$, $\kappa>0$, $a>0$. If
$\underline
\theta\!\in(0,2)$, then $\underline \ell(u) \sim
\underline \theta^{a-1} (-\log u)^{-a}$ as $u\to 0$. If
a L\'evy
process $X$ has $\nu^1_{a,\underline \theta}$ as a (symmetric) L\'evy
measure, then $r^*=+\infty$ and
\[
\forall\,
r,\,p\!\in(0,\underline \theta)\qquad  e_{N,r}(X,L^p_{T})
=
O((\log N)^{-{1}/{\underline\theta}}(\log\log N)^{-
a/2}).
\]
Such a rate improves the one provided
by~(\ref{BofRate0})

Let
$\nu^2_{a,\underline \theta}(dx) = \kappa|x|^{-\underline \theta
-1} (-\log |x|)^{a}\mbox{\bf 1}_{(0,c]}(|x|)\,dx$,
$\kappa,\,a>0$,
$0<c<1$,
 $\underline \theta\!\in(0,2)$.  Then $\underline \ell(u)
\sim \underline \theta^{-a-1} (-\log u)^{a}$ as $u\to 0$. Note  that
$\nu^2_{a,\underline
\theta}$ does not satisfy~(\ref{HypoLevy}).  If
a  L\'evy process $X$ has $\nu^2_{a,\underline \theta}$ as a
(symmetric) L\'evy measure, then
$r^*=+\infty$ and
\begin{eqnarray*}
\quad\hspace*{10pt}&&\forall\,
r,p\in(0,\underline \theta)\\
&&\qquad  e_{N,r}(X,L^p_{T})
=\cases{
O((\log N)^{-{1}/{\underline\theta}}(\log\log N)^{{a}/{(\underline
\theta -\eta)}}), \cr
\hspace*{147pt}\eta \in(0, \underline \theta), \mbox{ if $\underline\theta<1$,}\cr
O((\log N)^{-{1}/{\underline\theta}}(\log\log N)^{{2a}/{(\underline
\theta -\eta)}}),\cr
\hspace*{130pt}\eta  \in(0, \underline
\theta),\mbox{ if  $\underline \theta \in[1,2)$.}
}
\end{eqnarray*}

\textit{Hyperbolic
L\'evy motions} have been applied to option pricing in finance
(see~\cite{EBKE}). These processes are L\'evy processes whose
distribution $\mathbb{P}_{X_1}$
at time $1$ is a symmetric (centered)
hyperbolic distribution
\[
\mathbb{P}_{X_1}= C e^{-\delta
\sqrt{1+(x/\gamma)^2}}\,dx,\qquad \gamma,\delta>0.
\]
Hyperbolic L\'evy
processes are martingales with  no  Brownian component,   satisfying
$r^*=+\infty$. Their symmetric L\'evy measure has a Lebesgue density
that behaves like
$Cx^{-2}$ as $x\to 0$ [so that~(\ref{HypoLevy})
is satisfied with  $\underline \theta =1$]. Hence, one obtains, for
every $r,\, p\!\in (0,1)$,
\[
e_{N,r}(X,L^p_{T}) = O((\log
N)^{-1})
\]
and, for every $r\ge 1$ and every $p\in(0,r]$,
$e_{N,r}(X,L^p_{T}) = O((\log N)^{-1/r+\eta}),\break
\eta >0$.
\end{examples*}

The   proof of this theorem is divided into
several steps and is  deferred to Section~\ref{finalstep}. The reason
is that it
relies on the decomposition of $X$ as the sum of a
``bounded'' jump and a ``big'' jump L\'evy process. These are treated
successively in the  following
two sections.

%s4.1 ###
\subsection{L\'evy
processes with bounded jumps}\label{CompoundP}
In this section, we   consider  a  L\'evy process $X$
\textit{without Brownian component}
($\sigma=0$), with   jumps   bounded by a real
constant $c>0$.  In terms of the  L\'evy measure $\nu$ of $X$, this
means that
%
%e4.7 ###
\begin{equation}\label{Condic}
\nu([-c,c]^c)=0.
\end{equation}

Then, for every $\rho>0$ and every $t\ge 0$, $X_t\!\in
L^\rho(\mathbb{P})$,
that is, $r^*=+\infty$. In  Proposition~\ref{sous-theta}
below, we
establish   Theorem~\ref{Levyrate} in that
setting.

\begin{proposition}\label{sous-theta}
Let $(X_t)_{t\in[0,T]}$ be a
L\'evy process satisfying \textup{(\ref{Condic})} and\break $\underline \theta >0$.
Then
claims \textup{(a)}, \textup{(b)}, \textup{(c)} and \textup{(d)} in Theorem \textup{\ref{Levyrate}} hold
true with $r^*=\infty$.
\end{proposition}

\begin{pf}
The proof of this
proposition is decomposed into several steps.  We consider
$\underline
\theta$, as defined in Theorem~\ref{Upper1}. Note that, in
the present setting,\break
$\underline \theta = \inf\{\theta>0 \dvtx
\int|x|^\theta
\nu(dx)<+\infty\}$ and that $
\int|x|^\theta
\nu(dx)<+\infty$ for every $\theta >\underline
\theta$.  The
starting  point is to separate the ``small'' and the ``big''
jumps of $X$
in a nonhomogeneous way with respect to the
function
$s\mapsto s^{{1}/{\underline \theta}}$. We will
successively inspect the cases $\underline \theta \!\in(0,1)$ (or when $\underline
\theta =1$
holds as a minimum) and
$\underline
\theta\!\in[1,2]$.
\begin{step}[(\textit{Decomposition of $X$})]
 When $\underline \theta \!\in(0,1)$ or
$\underline \theta =1$ holds as a minimum, then
\[
\mathbb{E}\Bigg|\sum_{0<s\le
T}\Delta X_s\Bigg|  \le \mathbb{E}\sum_{0<s\le T}|\Delta X_s|  = T \int
|x|\nu(dx) <+\infty.
\]
Consequently, $X$ $\mathbb{P}$-a.s. has finite variation and we can decompose $X$ as
%e4.8 ###
\begin{equation}\label{Decompstart0}
X_t = \xi(t)+\sum_{0<s\le t} \Delta X_s,
\end{equation}
where $\xi(t)=at$ is a linear function.

Assume now  that $\underline \theta \!\in[1,2]$. We may decompose
$X$ as follows:
%e4.9 ###
\begin{eqnarray}\label{Decompstart}
X_t&=&\xi(t) +  X^{(\underline \theta)}_t + M^{(\underline \theta)}_t \quad \mbox{with} \nonumber\\
\xi(t)&:=&t\,\mathbb{E}(X_1), \\
X^{(\underline \theta)}_t&:=& \sum_{0<s\le t} \Delta
X_s\mathbf{1}_{\{|\Delta X_s|>s^{{1}/{\underline
\theta}}\}}-\int_0^t\,ds\!\int_{\{s^{{1}/{\underline \theta}}<
|x|\le c\}} x\nu(dx).\nonumber
\end{eqnarray}
Note that $X^{(\underline\theta)}$ has finite variations on $[0,T]$    since
\begin{eqnarray*}
&&\int_0^t ds\int_{\{s^{{1}/{\underline \theta}}< |x|\le
c\}}|x|\nu(dx)\\
&&\qquad= \int_{\{|x|\le
c\}}|x|(|x|^{\underline
\theta}\wedge t)\nu(dx)\le
\int_{\{|x|\le c\}} |x|^{1+\underline \theta} \nu(dx)<+\infty.
\end{eqnarray*}
Both $ X^{(\underline \theta)}$ and  $ M^{(\underline\theta)}$ are
martingales with (nonhomogeneous) independent
increments. Their   increasing predictable  ``bracket'' processes are given by
\[
 \bigl\langle X^{(\underline\theta)}\bigr\rangle_t=\int_0^t\!\!ds\!\int_{\{|x|>
s^{{1}/{\underline\theta}}\}}x^2\nu(dx)
\]
and
\[
\bigl\langle
M^{(\underline\theta)}\bigr\rangle_t=\int_0^t\,ds\int_{\{|x|\le
s^{{1}/{\underline\theta}}\}}x^2\nu(dx).
\]

From now on, we may consider  the (supremum process of the) L\'evy process
%
%e4.10 ###
\begin{equation}\label{tildeX}
\widetilde X_t := X_t-\xi(t),
\end{equation}
where $\xi$
is the linear function defined by~(\ref{Decompstart0})
and~(\ref{Decompstart}), respectively. Since the linear function $\xi$
lies in $L^p_{T}$,
it does not affect  the quantization rate, which  is invariant by translation.
\end{step}
\begin{step}[{[}Increment estimates in $L^\rho(\mathbb{P})${]}] In
this step, we   evaluate $\sup_{0\le s\le t} |\widetilde X_s|$  in
$L^\rho(\mathbb{P})$, $\rho\!\in(0,2]$. Throughout this step,  the $c$
comes from~(\ref{Condic}).
\end{step}
\begin{lema}\label{Lem2}
\textup{(a)}  Assume that $\underline\theta \in(0,1)$
or that $\underline \theta =1$ holds as a minimum. For every
$\rho\!\in(0,1]$
and
$t\in[0,T]$,
%
%e4.11 ###
\begin{eqnarray}\label{Ineqtheta1}
\qquad\mathbb{E} \biggl(\sup_{0\le s\le t} |\widetilde
X_s|^\rho\biggr)
&\le&  C_\rho\!
\biggl(\biggl(\int_0^t\int_{\{|x|\le
s^{{1}/{\underline \theta}}\}} x^2
\nu(dx)\biggr)^{{\rho}/{2}}\nonumber\\
&&\hspace{16pt}{}+\int_0^t\,ds\int_{\{s^{{1}/{\underline \theta}}< |x|\le c\}} |x|^\rho
\nu(dx)\\
&&\hspace{16pt}{}+  \sup_{0\le s\le t}\bigg|\int_0^s\,du
\int_{\{|x|\le   u^{{1}/{\underline \theta}}\}}x
\nu(dx)\bigg|^{\rho} \biggr).\nonumber
\end{eqnarray}

\textup{(b)} Assume that $\underline\theta \in[1,2]$.  For every
$\rho\!\in(0,2]$ and every $t\!\in[0,T]$,
%
%e4.12 ###
\begin{eqnarray}\label{Ineqtheta2}
\mathbb{E }\biggl(\sup_{0\le s\le t}|\widetilde   X_s|^\rho\biggr)
&\le&  C_\rho\biggl(\biggl(\int_0^t\,ds\int_{\{|x|\le
s^{{1}/{\underline \theta}}\}}x^2 \nu(dx)\biggr)^{{\rho}/{2}}\nonumber\\
&&\hspace{18pt}{}+\int_0^t\,ds\!\int_{\{s^{{1}/{\underline \theta}}< |x|\le c\}}
|x|^\rho \nu(dx)
\vid
&&\hspace{18pt}{} +\biggl(\int_0^t\,ds\! \int_{\{s^{{1}/{\underline
\theta}}< |x|\le c\}}|x|^{{\rho}/{2}}
\nu(dx)\biggr)^{2}\biggr)\nonumber\\
&&{}+ \sup_{0\le
s\le t}\bigg|\int_0^s\,
du\int_{\{u^{{1}/{\underline \theta}}< |x|\le c\}}
x\nu(dx)\bigg|^\rho.\nonumber
\end{eqnarray}
\end{lema}

\begin{pf}
(a) $\widetilde X$ is a pure jump
process (with finite variations). Using   $\rho\!\in(0,1]$
and Doob's inequality, we obtain
\begin{eqnarray*}
&&\mathbb{E} \sup_{0\le s\le t} |\widetilde X_s|^\rho\\
&&\qquad\le  \mathbb{E} \sup_{0\le
s\le t} \Bigg|\sum_{0\le u\le s} \Delta
X_u\mbox{\bf 1}_{\{|\Delta X_u|\le u^{{1}/{\underline \theta}} \}
}\Bigg|^\rho+ \mathbb{E} \!\sup_{0\le s\le t} \Bigg|\sum_{0\le u\le s} \Delta
X_u\mathbf{1}_{\{|\Delta X_u|> u^{{1}/{\underline \theta}} \} }\Bigg|^\rho\\
&&\qquad\le \Biggl(\mathbb{E} \!\sup_{0\le s\le t} \Biggl(\sum_{0\le u\le s} \Delta
X_u\mathbf{1}_{\{|\Delta X_u|\le u^{{1}/{\underline \theta}} \}
}\Biggr)^2\Biggr)^{{\rho}/{2}}+ \mathbb{E} \!\sum_{0<s\le t}
|\Delta X_s|^{\rho}\mathbf{1}_{\{|\Delta X_s|>
s^{{1}/{\underline \theta}} \} }\\
&&\qquad\le C_\rho\Biggl(\Biggl(\mathbb{E} \!\sup_{0\le s\le t} \Biggl(\sum_{0\le u\le s} \Delta
X_u\mathbf{1}_{\{|\Delta X_u|\le u^{{1}/{\underline \theta}} \}
}-\int_0^s\,du\int_{\{|x|\le   u^{{1}/{\underline
\theta}}\}}\!\!x\nu(dx)\Biggr)^2\Biggr)^{{\rho}/{2}}\\
&&\qquad\quad\hspace{24pt}{}+ \sup_{0\le s\le t}\bigg|\int_0^s\,du\int_{\{|x|\le
u^{{1}/{\underline
\theta}}\}}\!\!x\nu(dx)\bigg|^{\rho}+\int_0^t\,ds\int_{\{s^{{1}/{\underline
\theta}}< |x|\le c\}}\!\!|x|^\rho\nu(dx)\Biggr)\\
&&\qquad\le C_\rho\Biggl(\Biggl(\int_0^t\,ds\!\!\int_{\{|x|\le
s^{{1}/{\underline
\theta}}\}}\!\!x^2\nu(dx)\Biggr)^{{\rho}/{2}}\\
&&\qquad\quad\hspace{24pt}{}+ \sup_{0\le s\le
t}\bigg|\int_0^s\,du\int_{\{|x|\le   u^{{1}/{\underline
\theta}}\}}\!\!x\nu(dx)\bigg|^{\rho}+\int_0^t\,ds\!\!\int_{\{s^{{1}/{\underline
\theta}}< |x|\le c\}}\!\!|x|^\rho\nu(dx)\Biggr).
\end{eqnarray*}

(b)  It follows from   Doob's inequality (and $0<\rho/2\le 1$) that
\begin{eqnarray*}
\mathbb{E}\biggl(\sup_{0\le s\le t}\big|M^{(\underline\theta)}_s\big|^\rho  \biggr)&\le& \biggl[\mathbb{E}\sup_{0\le
s\le t} \bigl(M^{(\underline \theta)}_s\bigr)^2\biggr]^{{
\rho}/{2}}\le \biggl(4\int_0^t\,ds\! \int_{\{|x|\le
s^{{1}/{\underline \theta}}\}}x^2 \nu(dx)\biggr)^{{
\rho}/{2}}.
\end{eqnarray*}

On the other hand, since $\rho \!\in(0,2]$, we have
\begin{eqnarray*}
&&\sup_{0\le s\le t} \big|X_s^{(\underline \theta)} \big|^\rho\\
&&\qquad\le C_\rho\Biggl( \Biggl(\sum_{0<s\le t} |\Delta
X_s|^{{\rho}/{2}}\mathbf{1}_{\{|\Delta
X_s|> s^{{1}/{\underline \theta}} \}}\Biggr)^2\\
&&\qquad\hspace{18pt}{}+ \sup_{0\le s\le
t}\bigg|\int_0^s\,du\!\int_{\{u^{{1}/{\underline \theta}}<
|x|\le c\}}
x\nu(dx)\bigg|^\rho \Biggr)\\
&&\qquad\le C_\rho\Biggl(\Biggl(\sum_{0<s\le t} |\Delta
X_s|^{{\rho}/{2}}\mathbf{1}_{\{|\Delta X_s|>
s^{{1}/{\underline \theta}} \}}-\int_0^t\,ds
\int_{\{|x|>  s^{{1}/{\underline \theta}}\}}|x
|^{{\rho}/{2}}\nu(dx) \Biggr)^2\\
&&\qquad\hspace{30 pt{}}{}+\biggl(\int_0^t\,ds\!
\int_{\{|x|>  s^{\frac{1}{\underline \theta}}\}}|x
|^{{\rho}/{2}}\nu(dx) \biggr)^2\\
&& \qquad\hspace{140pt}{}+   \sup_{0<  s\le
t}\bigg|\int_0^s\,du\int_{\{u^{{1}/{\underline \theta}}<
|x|\le c\}}
x\nu(dx)\bigg|^\rho\Biggl).
\end{eqnarray*}
Hence, again using Doob's inequality,
\begin{eqnarray*}
&&\mathbb{E} \!\sup_{0\le s\le t}
\big|X^{(\underline\theta)}_s\big|^{\,\rho}\\
&&\qquad\le
C_\rho \biggl(\int_0^t\,ds
\int_{\{s^{{1}/{\underline \theta}}< |x|\le c\}}|x
|^{\rho}\nu(dx) +  \biggl(\int_0^t\,ds
\int_{\{s^{{1}/{\underline \theta}}< |x|\le c\}}|x
|^{{\rho}/{2}}\nu(dx)\biggr)^2\\
&&\hspace{178pt}{}+\sup_{0<  s\le t}\biggl|\int_0^s\,du\int_{\{u^{{1}/{\underline
\theta}}< |x|\le c\}}
x\nu(dx)\biggr|^\rho\biggr).
\end{eqnarray*}
\hspace{337pt}\ \qed\noqed
\end{pf}\noqed
\begin{lema}[(First extended Millar's lemma)]\label{Quantifrate1}
 \textup{(a)}
Assume that $\underline\theta \!\in (0,2]\setminus\{1\}$. If   the L\'evy
measure
satisfies assumption~(\ref{HypoLevy})  then
%e4.13 ###
\begin{equation}\label{BofRatemart}
\forall\, \rho\in(0,\underline \theta ),  \; \forall\, t\in[0,
T] \qquad  \mathbb{E} \sup_{0\le s \le t} |\widetilde X_s|^\rho
\le    C_\rho t^{{\rho}/{\underline
\theta}}.
\end{equation}

{\smallskipamount=0pt
\begin{longlist}[(b)]
\item[(b)] Assume that $\underline\theta \!\in (0,2)\setminus\{1\}$
and  that the function $u\mapsto \underline \nu(u)$ has regular
variation with index $-b$ at $0$. Then $b=\underline
\theta$ and, for every $ \rho\!\in(0,\underline \theta )$, there
exists $T_\rho\!\in (0,T]$ such that
%e4.14 ###
\begin{equation}\label{BofRate2}
\forall\, t\!\in[0,T_\rho ] \qquad
\mathbb{E} \sup_{0\le s \le t} |\widetilde X_s|^\rho  \le  C_\rho
(t^{{1}/{\underline\theta}} \Lambda_\rho(t))^\rho.
\end{equation}

\item [(c)] When $\underline \theta =1$, the above upper bounds
still hold,   provided $\nu$ is symmetric   or  $\nu(|x|)<+\infty$.
\end{longlist}
}
\end{lema}

\begin{pf}
(a) We need only to investigate all the
integrals appearing in the right-hand side of
inequalities~(\ref{Ineqtheta1}) and~(\ref{Ineqtheta2})  in
Lemma~\ref{Lem2}.  Let $\rho \!\in (0,\underline \theta)$ and
$t\!\in[0,c^{\,\underline
\theta}\wedge T]$. Then, if $\underline \theta \!\in(0,2)$,
\begin{eqnarray*}
\int_0^t  \!ds\int_{\{0<|x|\le s^{{1}/{\underline
\theta}}\}}x^2\nu(dx)&\le& C \int_0^t\,ds\int_{\{0<|x|\le
s^{{1}/{\underline
\theta}}\}}|x|^{1-\underline \theta}\,dx\\
&\le& C \int_0^t s^{{2}/{\underline \theta}-1}\,ds = C
t^{{2}/{\underline\theta}},
\end{eqnarray*}
where the real constant $C$ comes from~(\ref{HypoLevy}). If
$\underline \theta =2$, then
\[
\int_0^t  \,ds\int_{\{0<|x|\le s^{{1}/{\underline
\theta}}\}}x^2\nu(dx)\le  \int_{[-c,c]}x^2\nu(dx)t
= \int_{[-c,c]} x^2\nu(dx)t^{{2}/{\underline \theta}}.
\]
Then, for every
$t\!\in[0,c^{\,\underline \theta}\wedge T]$,
\begin{eqnarray*}
\int_0^t  \,ds\int_{\{s^{{1}/{\underline
\theta}}<|x|\le c\}}|x|^\rho\nu(dx)&\le& C
\int_0^t\,ds\int_{\{s^{{1}/{\underline
\theta}}<|x|\le c\}}|x|^{\rho-\underline
\theta-1}\,dx\\
&\le& {C}/{\underline \theta -\rho} \int_0^t
s^{{\rho}/{\underline \theta}-1} \,ds= C t^{{\rho}/{\underline
\theta}}.
\end{eqnarray*}
When $\underline \theta \!\in(0,1)$, we have
\begin{eqnarray*}
\sup_{0\le s\le t}\bigg|\int_0^s \!du\int_{\{|x|
\le u^{{1}/{\underline
\theta}}\}}x\nu(dx)\bigg|
&\le& \int_0^t
ds\int_{\{|x|\le s^{{1}/{\underline
\theta}}\}}|x|\nu(dx)\\
&\le& C \int_0^t \frac{ s^{{1}/{\underline
\theta}-1}}{1-\underline \theta}\,ds = \frac{C}{{1-\underline \theta}}
t^{{{1}/{\underline \theta}}}.
\end{eqnarray*}
When $\underline \theta= 1$ and $\int |x|\nu(dx)<+\infty$, this
term is trivially upper bounded by $t\int |x|\nu(dx)$. It  is $0$
when $\nu$
is symmetric.  Similarly, when
$\underline
\theta
\!\in(1,2]$, for every
$t\!\in[0,c^{\,\underline \theta}\wedge T]$, we have
\begin{eqnarray*}
\sup_{0\le s\le t}\bigg|\int_0^s \,du\int_{\{u^{{1}/{\underline
\theta}}<|x|\le c\}}x\nu(dx)\bigg|&\le& \int_0^t
\,ds\!\int_{\{|x|> s^{{1}/{\underline \theta}}\}}|x|\nu(dx)\\
&\le& C
\int_0^t
\frac{ s^{{1}/{\underline \theta}-1}}{\underline \theta-1}\,ds =
\frac{C}{{\underline \theta}-1} t^{{{1}/{\underline
\theta}}}
\end{eqnarray*}
and
\begin{eqnarray*}
\int_0^t  \,ds\int_{\{s^{{1}/{\underline
\theta}}<|x|\le c\}}|x|^{{\rho}/{2}}\nu(dx)
&\le& C
\int_0^t\,ds\int_{\{s^{{1}/{\underline
\theta}}<|x|
\le c\}}|x|^{{\rho}/{2}-\underline
\theta-1}\,dx\\
&\le& \frac{C}{\underline \theta -{\rho}/{2}} \int_0^t
s^{{\rho}/({2\underline \theta})-1} \,ds= C t^{{\rho}/({2\underline
\theta})}.
\end{eqnarray*}
It can be derived from~(\ref{Ineqtheta1}) and~(\ref{Ineqtheta2}) that there
exists  a positive real constant  $C_\rho$ such that
\begin{eqnarray*}
\forall\,  t\in[0,c^{\,\underline \theta}\wedge T]\qquad \mathbb{E}
\sup_{0\le s \le t} |\widetilde X_s|^\rho &\le &
C_\rho t^{{\rho}/{\underline
\theta}}.
\end{eqnarray*}
This inequality  holds for every $t\in [0,T]$ simply by adjusting
the constant $C_\rho$.

(b) The fact that $b= \underline \theta$  was first
established in~\cite{BLGE}. We provide below a short proof, leading
to our main result, for the reader's convenience. It follows from
Theorem 1.4.1 in~\cite{BIGOTE} that
$\underline
\nu(u)= u^{-b}\ell(u)$ where
$\ell$ is a (nonnegative) slowly varying function.  Consequently,
one clearly has  that, for every $\rho>0$ and every $u>0$,
\[
  u^{\rho-b} \ell(u) \le \int_{\{|x|>u\}}|x|^\rho \nu(dx).
\]
Now, the left-hand side of the above inequality goes to infinity as
$u\to 0$ provided $\rho < b$ since $\ell$ has slow
variations (see Proposition~1.3.6 in~\cite{BIGOTE}). Consequently,
$\rho\le \underline \theta$. Letting $\theta$ go to $b$ implies that
$b\le \underline\theta$.

We will make use  of the following easy  identity which follows from
the very definition of
$\underline\nu$: for every nonnegative Borel function
$f\dvtx\mathbb{R}_+\to\mathbb{R}$,
%e4.15 ###
\begin{equation}\label{identitefnu}
\int_{\mathbb{R}} f(|x|)\nu(dx) = - \int_{\mathbb{R}_+} f(x)\,d\underline \nu(x).
\end{equation}

In particular, for every $x\!\in (0,c]$ and every $a>0$,
\[
\int_{\{|u|\ge x\}}|u|^a  \nu(du) =- \int_x^c u^a d\underline\nu(u).
\]
Assume that  $b<\underline \theta$. It then follows from Theorem 1.6.4
in~\cite{BIGOTE} that for every $a\!\in(b,\underline \theta)$,
\[
\int_x^c u^a\, d\underline\nu(u) \sim \frac{b}{b-a}\, x^a \,\underline
\nu(x)=\frac{b}{b-a}\, x^{a-b}  \ell(x) \to 0\qquad\mbox{as }
x\to 0,
\]
since $\ell$ is slowly varying. This contradicts $\int|u|^a  \nu(du)
=+\infty$. Consequently,\break $b=\underline
\theta$.

Now,   Theorem~1.6.5 in~\cite{BIGOTE} implies that for any
$a>\underline \theta$
\[
\int_{\{|u|\le x\}} |u|^a \nu(du)  =- \int_{(0,x]} u^a
\,d\underline\nu(u)\sim \frac{\underline \theta}{a-\underline \theta}
x^a \underline \nu(x)\qquad \mbox{as }x\to 0.
\]
Since $\underline \theta \neq 2$, this yields
\[
\int_{\{|x|\le s^{{1}/{\underline \theta}}\}} x^{2}\nu(dx) \sim
\frac{\underline \theta}{2-\underline \theta}
\,s^{{2}/{\underline \theta}}\,\underline \nu(s^{{1}/{\underline
\theta}})\qquad \mbox{as }s\to 0,
\]
which, in turn, implies that
\[
\int_0^t ds\int_{\{|x|\le s^{{1}/{\underline \theta}}\}}
x^{2}\nu(dx) \sim \frac{\underline \theta}{2-\underline \theta}
\int_0^ts^{{2}/{\underline \theta}}\,\underline
\nu(s^{{1}/{\underline
\theta}})\,ds\qquad \mbox{as }t\to 0.
\]
The function $s\mapsto\underline \nu(s^{{1}/{\underline
\theta}})$ has regular variation (at $0$) with index $-1$, hence
Theorem~1.6.1 in~\cite{BIGOTE} implies that
\[
\int_0^t ds\int_{\{|x|\le s^{{1}/{\underline \theta}}\}}
x^{2}\nu(dx) \sim C_{\underline \theta}\,
t^{{2}/{\underline \theta}+1} \underline \nu(t^{{1}/{\underline
\theta}})\qquad
\mbox{as } t\to 0.
\]
Finally,
%e4.16 ###
\begin{equation}\label{equivat0}
\qquad\biggl(\int_0^t ds\int_{\{|x|\le s^{{1}/{\underline \theta}}\}}
x^{2}\nu(dx)\biggr)^{{\rho}/{2}} \sim
C_{\rho,\underline \theta}(t^{{1}/{\underline
\theta}}  (\underline
\ell(t))^{{1}/{2}})^\rho  \qquad
\mbox{as } t\to 0.
\end{equation}
When $\underline \theta \!\in(0,1)$ and
$\rho\!\in(0,\underline\theta)$,   the same approach  leads to
\begin{eqnarray*}
&&\sup_{0< s\le t}\bigg|\int_0^s du\int_{\{|x|\le
u^{{1}/{\underline \theta}}\}} x\nu(dx)\bigg| \\
&&\qquad\le  \int_0^t
ds\int_{\{|x|\le s^{{1}/{\underline
\theta}}\}} |x|\nu(dx)  \sim C_{\underline \theta}\,
  t^{{1}/{\underline \theta}} \underline \ell(t) \qquad
\mbox{as } t\to 0.
\end{eqnarray*}
It then follows from   Theorem~1.6.4 in~\cite{BIGOTE}
that, for every
$\rho\!\in (0,\underline\theta)$,
\begin{eqnarray*}
\int_{\{s^{{1}/{\underline \theta}}\le |x|\le c \}} |x|^\rho
\nu(dx)&=&-\int_{s^{{1}/{\underline \theta}}}^c x^\rho\,d\underline \nu(x)\\
&\sim& \frac{\underline \theta}{\underline
\theta-\rho } s^{{\rho}/{\underline \theta}} \underline
\nu(s^{{1}/{\underline
\theta}})\qquad \mbox{as }s\to 0
\end{eqnarray*}
so that
\begin{eqnarray*}
\int_0^t\int_{\{s^{{1}/{\underline
\theta}}\le |x|\le c \}} |x|^{\rho}  \nu (dx)  &\sim&  {\underline
\theta}/{\underline
\theta-\rho }\int_0^ts^{{\rho}/{\underline
\theta}}
\underline \nu(s^{{1}/{\underline \theta}})\,ds\\
&\sim&
C_{\rho,\underline \theta}\, (t^{{1}/{\underline \theta}}
(\underline
\ell(t) )^{{1}/{\rho}})^\rho\qquad \mbox{as } t\to 0.
\end{eqnarray*}

Similarly (by formally setting $\rho=1$ in the former equation)  we can shown
that if $\underline \theta\!\in(1,2]$, then
%e4.17 ###
\begin{eqnarray}\label{equivat1}
\qquad\sup_{0< s\le t}\bigg|\int_0^s \,du\int_{\{u^{{1}/{\underline
\theta}}\le |x|\le c \}}\!\! x\nu(dx)\bigg| &\le& \int_0^t
ds\int_{\{s^{{1}/{\underline \theta}}\le |x|\le c \}}
|x|\nu(dx)
\vid
&\sim& C_{\underline \theta}\,
  t^{{1}/{\underline \theta}} \underline \ell(t) \qquad
\mbox{as }t\to 0.\nonumber
\end{eqnarray}
Finally, we similarly shown, for   the last term
in~(\ref{Ineqtheta2}), that when
$\rho\!\in(0,\underline
\theta)$,
\[
\biggl(\int_0^t\int_{\{s^{{1}/{\underline \theta}}\le |x|\le c \}}
|x|^{{\rho}/{2}}  \nu (dx)\biggr)^2 \sim   C_{\rho,\underline
\theta}\,(t^{{1}/{\underline
\theta}} (\underline \ell(t))^{{2}/{\rho}})^\rho \qquad
\mbox{as }t\to 0.
\]

Substituting these estimates into (\ref{Ineqtheta1}) and~(\ref{Ineqtheta2})
and noting that, by Young's inequality,
\[
\underline \ell(t) \le  C_\rho\bigl((\underline
\ell(t))^{{1}/{2}}+(\underline \ell(t))^{{1}/{\rho}}\mathbf{1}_{\{\rho\le
1\}}+(\underline
\ell(t))^{{2}/{\rho}}\mathbf{1}_{\{1<\rho\le 2\}}\bigr),
\]
we finally obtain that $\widetilde X$ satisfies
the assumption
$(L_{\varphi,\rho})$ with the announced function $\varphi_{\rho}$.

(c) When $\nu$ is symmetric (and   $\underline \theta
\!\in(1,2]$), for every $s\!\in [0,T]$,
\[
\int_0^s \!du\!\int_{\{u^{{1}/{\underline \theta}}\le |x|\le c
\}} x\nu(dx)=0
\]
so that the condition $\underline \theta \neq 1$ induced
by~(\ref{equivat1}) is no longer necessary. Similarly, when
$\underline \theta \!\in(0,1]$,
\[
\int_0^s \!du\!\int_{\{|x|\le u^{{1}/{\underline \theta}}\}}
x\nu(dx)= 0.
\]\upqed
\end{pf}

\begin{step}[(Higher moments and completion of the proof)] Claims
(a), when $\underline \theta$ holds as a minimum, and (c), when
$r<  2$, straightforwardly follow from
Millar's inequality~(\ref{Millar})  by applying Theorem~\ref{Upper1}
to the function $\varphi(u) = u^{{1}/{\underline \theta}}$ with
$\rho
=\underline \theta$ for claim
(a) and $\varphi(u) = u^{{1}/{\rho}}$ with $\rho\!\in(r,2]$ for
claim~(c).

Claim~(a), when assumption~(\ref{HypoLevy})  is
fulfilled, follows from Lemma~\ref{Quantifrate1}(a) and
Theorem~\ref{Upper1} applied with  the function $\varphi(u) =
u^{{1}/{\underline
\theta}}$. Finally, claim~(b) follows from
Lemma~\ref{Quantifrate1}(b) and Theorem~\ref{Upper1}.

Claim~(d) follows from Lemma~\ref{Quantifrate1}(c) and Theorem~\ref{Upper1}.
At this stage, it remains to prove claim~(c) when $r\ge  2$. This
follows (when $r>2$) from the extension of Millar's upper bound
established in the lemma
below.
\end{step}

\begin{lema}[(Second extended Millar's lemma)]
Let
$(X_t)_{t\in[0,T]}$ be a L\'evy process without Brownian part such
that $\nu([-c,c]^c)=0$. For every $\rho\ge 2$,
  there exists  a real constant $C_{\rho,T}>0$ such
that
\[
\forall\, t\in[0,T] \qquad \mathbb{E}\biggl(\sup_{0\le s\le t}|X_s|^\rho\biggr)\le
C_{\rho,T} \,t.
\]
\end{lema}

\begin{pf}
We again consider $\widetilde X_t= X_t- t
\,\mathbb{E} X_1$, which is a martingale L\'evy process.  Let $k_\rho :=
\max\{l\dvtx 2^l <\rho\}$. For every
$k=1,\ldots,k_\rho$, we define the martingales
\[
N^{(k)}_t := \sum_{0<s\le t} |\Delta X_s|^{2^k} -t \int |x|^{2^k} \nu(dx).
\]
The key technique of the proof is to apply the BDG inequality in cascade. It
follows from the BDG inequality that
\begin{eqnarray*}
\mathbb{E} \sup_{0\le s\le t} |\widetilde X_s|^\rho  &\le &C_\rho
\mathbb{E}\Biggl(\sum_{0<s\le t}(\Delta X_s)^2\Biggr)^{\rho/2}\\
&\le& C_\rho\biggl( \mathbb{E}\bigl(N^{(1)}_t\bigr)^{\rho/2}+ \biggl(t\int
x^2\nu(dx)\biggr)^{\rho/2}\biggr).
\end{eqnarray*}

Now, for every $k\!\in\{1,\ldots,k_\rho-1\}$, still using the BDG
inequality yields
\begin{eqnarray*}
\mathbb{E}\bigl(N^{(k)}_t\bigr)^{\rho/2^k} &\le& C_{\rho,k}
\mathbb{E}\Biggl(\sum_{0<s\le t}|\Delta
X_s|^{2^{k+1}}\Biggr)^{\rho/2^{k+1}}\\
&\le & C_{\rho,k}\biggl( \mathbb{E}\bigl(N^{(k+1)}_t\bigr)^{\rho/2^{k+1}}+ \biggl(t\int
|x|^{2^{k+1}}\nu(dx)\biggr)^{\rho/2^{k+1}}\biggr).
\end{eqnarray*}
Finally, we obtain
\begin{eqnarray*}
\mathbb{E} \sup_{0\le s\le t} |\widetilde X_s|^\rho  &\le &C_\rho
\Biggl(\sum_{k=1}^{k_\rho} \biggl(t\int
|x|^{2^{k}}\nu(dx) \biggr)^{\rho/2^k}+\mathbb{E}  \Biggl(\sum_{0<s\le
t}|\Delta X_s|^{2^{k_\rho+1}}\Biggr)^{\rho/2^{k_\rho+1}}\Biggr)\\
&\le& C_\rho  \Biggl(\sum_{k=1}^{k_\rho} \biggl(t\int
|x|^{2^{k}}\nu(dx) \biggr)^{\rho/2^k}+ \mathbb{E} \sum_{0<s\le t}|\Delta
X_s|^{\rho} \Biggr)\\
&=& C_\rho  \Biggl(\sum_{k=1}^{k_\rho} \biggl(t\int
|x|^{2^{k}}\nu(dx) \biggr)^{\rho/2^k}+ t\int|x|^\rho\nu(dx) \Biggr)
\end{eqnarray*}
since $\rho/2^{k_\rho+1}\le 1$. The conclusion follows from the fact
that $t^{\rho/2^k} = o(t)$.\quad\qed
\end{pf}\noqed
\end{pf}\noqed

%s4.2 ###
\subsection{Compound  Poisson process}\label{CompP2}
In this section,
we consider a compound Poisson process $(X_t)_t$
defined by
\[
X_t :=  \sum_{n\ge 1}U_n\mathbf{1}_{\{S_n\le \lambda T\}},\qquad t\ge0,
\]
where  $S_n =Z_1+\cdots+Z_n$, $(Z_n)_{n\ge 1}$ is an i.i.d. sequence
of $\mathcal{E}\!xp(1)$-distributed random
variables, $(U_n)_{n\ge 1}$ is an  i.i.d. sequence of
random variables, independent of $(Z_n)_{n\ge 1}$ with $ U_1\!\in
L^\rho$, $\rho>0$ and  $\lambda >0$
  is the the jump intensity. For convenience, we also introduce the
underlying standard Poisson process $(K_t)_{t\ge 0}$ defined by
\[
K_t:=  \sum_{n\ge 1}\mathbf{1}_{\{S_n\le \lambda T\}},\qquad t\ge 0,
\]
so that (with the convention that $\sum_{\varnothing} =0$)
%e4.18 ###
\begin{equation}\label{PoissonComp}
X_t =\sum_{k=1}^{K_t}  U_k.
\end{equation}

\begin{proposition}\label{Poisson}
Let $X$ be a compound Poisson process.
Then, for every $p,r\!\in(0, r^*)$, $p\le r$,
%
%e4.19 ###
\begin{eqnarray}\label{ratePoisson}
\qquad&&\forall\, \varepsilon>0
\vid
&&\qquad e_{N,r}(X,L^p_{T})=
O\biggl(\exp{\biggl(-\frac{1}{\sqrt{r(p+1+\varepsilon)}}\sqrt{\log(N)\log_2(N)}\biggr)}\biggr).\nonumber
\end{eqnarray}
Furthermore, when $X$ is a standard Poisson process, we can
replace $p+1+\varepsilon$ by
$p+\varepsilon$ in~(\ref{ratePoisson}).
\end{proposition}

\begin{remark*}
Note that~(\ref{ratePoisson}) implies that
\[
\forall\, a>0\qquad e_{N,r}(X,L^p_{T}) = o((\log N)^{-a}).
\]

In fact, the rate obtained in the above
proposition holds provided $X$ has the
form~(\ref{PoissonComp}), where $(Z_n)$ is as above and  $(U_n)$
is $L^r(\mathbb{P})$-bounded for every $r<r^*$, independent of $(Z_n)_{n\ge 1}$.
\end{remark*}

\begin{pf*}{Proof of Proposition \protect\ref{Poisson}}
We divide the proof into two steps, one devoted to
the standard Poisson process, the
other to the general case. We will assume that $r^*> 1$
throughout the proof so that, as
was already emphasized in  the proof of Theorem~\ref{Upper1},  we may
assume without loss of generality that
$r,\,p\!\in (0,r^*)\cap [1,+\infty)$. The case $r^*\le 1$ is left to
the reader, but can be treated by
replacing the ``triangular'' Minkowski inequality  by the
pseudo-triangular inequalities
$
|f+g|_{L^p_{T}}^p\le |f|_{L^p_{T}}^p+|g|_{L^p_{T}}^p
$
and
$\|U+V\|_{r} ^r \le \|U\|_{r} ^r+\|V\|_{r} ^r$.
\setcounter{step}{0}
\begin{step}[(Standard case)] One quantizes the standard
Poisson $K$ in a very natural way by setting
\[
\widehat K_t := \sum_{n\ge 1} \mathbf{1}_{\{\widehat S_n\le \lambda
t\}},\qquad t\ge 0,
\]
with
\[
\widehat S_n:= \widehat{S_n}^{\alpha_n},
\]
where $\alpha_n = \alpha'_n\cup\{\lambda T\}$, $\alpha'_n$ is an
$L^{r'}$-optimal  $(N_n-1)$-quantization
of
$S_n^{tr}:=S_n\mathbf{1}_{\{S_n\le  \lambda T\}}$ and  $r'= \frac
rp$.  Furthermore, we assume that the
sequence
$(N_n)$ is nonincreasing and satisfies
$\prod_n N_n
\le N$ (so that $N_n=1$ for large enough $n$). Then, for every
$p\ge 1$, it follows from the (extended) Minkowski  inequality that
\[
|K- \widehat K|_{L^p_{T}}\le \sum_{n\ge 1}\bigl|\mathbf{1}_{\{S_n\le
\lambda \bolds{\cdot}\}}-\mathbf{1}_{\{\widehat S_n\le \lambda \bolds{\cdot}\}}\bigr|_{L^p_{T}}.
\]
Now,
\begin{eqnarray*}
\big|\mathbf{1}_{\{S_n\le \lambda \cdot\}}-\mathbf{1}_{\{\widehat S_n\le
\lambda \cdot\}}\big|^p_{L^p_{T}}
&=& \int_0^T\big|\mathbf{1}_{\{S_n\le \lambda
t\}}-\mathbf{1}_{\{\widehat S_n\le \lambda t\}}\big|^p\,dt\\
&=&   \frac {1} {\lambda} \,|S_n \wedge(\lambda T)- \widehat
S_n\wedge(\lambda T)|= \frac 1 \lambda \,|S_n \wedge(\lambda T)- \widehat
S_n|.
\end{eqnarray*}

Now,   $\{S_n> \lambda T\}\subset\{\widehat S_n= \lambda
T\}$ since $\max\alpha_n = \lambda T$. On
the other hand, $S_n=S_n^{tr}$ on $\{S_n\le  \lambda T\}$ so that
\[
|S_n \wedge(\lambda T)- \widehat S_n|= |S_n \wedge(\lambda T)-
\widehat S_n|\mathbf{1}_{\{S_n \le \lambda
T\}}= |S^{tr}_n  - \widehat{S^{tr}_n}|\mathbf{1}_{\{S_n \le \lambda
T\}}\le   |S^{tr}_n  - \widehat{S^{tr}_n}|.
\]
Also, note that  when $N_n=1$,  $\widehat S_n = \lambda T$ so that
$|S_n \wedge(\lambda T)- \widehat
S_n|=(\lambda
T-S_n)_+$. Consequently, for every $r\ge 1$,
\begin{eqnarray*}
\big\||K-\widehat K|_{L^p_{T}}\big\|_{r}&\le & \sum_{n\ge
1}\,\big\|\big|\mathbf{1}_{\{S_n\le \lambda \cdot\}}-\mathbf{1}_{\{\widehat
S_n\le \lambda
\cdot\}}\big|_{L^p_{T}}\big\|_{r}\\
&\le & \dfrac{1}{\lambda^{1/p}}\sum_{n\ge 1}\|S_n
\wedge(\lambda T)- \widehat S_n\|_{{r'}}^{
1/p}\\
&\le & \frac {1}{\lambda^{1/p}}\Biggl(\sum_{n, N_n\ge 2}\,\|S^{tr}_n  -
\widehat{S^{tr}_n}^{\alpha_n}\|^{1/p}_{{r'}}+\sum_{n,
N_n=1}\,\|(\lambda T-S_n)_+\|^{
1/p}_{{r'}}\Biggr)\\
 &\le& \frac {1}{\lambda^{1/p}}\Biggl(\sum_{n, N_n\ge
2}\,\|S^{tr}_n  -\widehat{S^{tr}_n}^{\alpha'_n}\|^{1/p}_{{r'}}+\sum_{n,
N_n=1}\,\|(\lambda T-S_n)_+\|^{1/p}_{{r'}}\Biggr).
\end{eqnarray*}
The extended Pierce lemma~(Lemma~\ref{UPierce}) yields  that, for every
$n\ge 1$ such that $N_n\ge 2$ and  for every $\delta>0$,
\begin{eqnarray*}
\|S^{tr}_n  -\widehat{S^{tr}_n}^{\alpha'_n}\|_{{r'}}
&\le & \|S_n^{tr}\|_{{r'+\delta/p}} C_{r,p,\delta} \,|N_n-1|^{-1}\\
&\le& 2 \big\|S_n\mathbf{1}_{\{S_n\le  \lambda T\}}\big\|_{({r+\delta})/{p}}
C_{r,p,\delta} \,N_n^{-1}.
\end{eqnarray*}

Set $\mu :=  r'+\delta/p =\frac{r +\delta}{p}$ so that $\mu p=r+\delta$.
We then have
%e4.20 ###
\begin{eqnarray}\label{quantN}
\big\||K-\widehat K|_{L^p_{T}}\big\|_{r} &\le&
C_{p,r,\delta}\frac {1}{\lambda^{1/p}}\Biggl(\sum_{n, N_n\ge
2}\big\| S_n\mathbf{1}_{\{ S_n\le \lambda
T\}}\big\|^{1/p}_{{\mu}}\frac{1}{N^{ 1/p}_n}\nonumber\\
&&\hspace{64pt}{}+\sum_{n,
N_n=1}\,\|(\lambda T-S_n)_+\|^{ 1/p}_{{\mu}}\Biggr)\\
&\le& C_{ p,r,\delta}T^{1/p}\Biggl(\sum_{n\ge
1}\,\bigl(\mathbb{P}( S_n\le \lambda T)\bigr)^{{1}/{(\mu p)}}
\frac{1}{N^{ 1/p}_n}\Biggr).\nonumber
\end{eqnarray}
Now, standard computations show that
\[
\mathbb{P}(\{ S_n\le \lambda T\})  = \frac{(\lambda
T)^{n}}{(n-1)!}\int_0^1u^{n-1}e^{-\lambda T u}\,du\le \frac{(\lambda
T)^{n}}{n!}.
\]
Hence, setting $A = (\lambda T)^{{1}/({\mu p})}$ yields
\[
\bigl(\mathbb{P}( S_n\le \lambda T)\bigr)^{{1}/({\mu p})}\le
\frac{(\lambda T)^{{n}/({\mu p})}}{(n!)^{{1}/({\mu p})} }\le
\frac{A^n}{(n!)^{{1}/({\mu p})}}.
\]

For every $x\ge 0$, let  $a(x):= \frac{A^x}{\Gamma(x+1)^{{1}/({\mu
p})}}$. This function reaches a unique maximum at
some
$x_0\ge 0$ and then decreases to $0$ as
$x\to \infty$. We modify the function $a$ by setting $  a_0(x) :=
a(x)\vee a(x_0)$ so that the function $a_0$ becomes
nonincreasing and $\log$-concave since $\Gamma$ is $\log$-convex. Now, let
\[
a_n :=   a_0(n),\qquad n\ge 1.
\]

Finally, the quantization problem~(\ref{quantN}) for the standard
Poisson $K$  is ``upper bounded'' by the following  optimal integral
``bit allocation'' problem:
%e4.21 ###
\begin{equation}\label{bitalloc}
\min\Biggl\{\sum_{n\ge 1}\ \frac{a_n}{N^{1/p}_n},\; N_n\ge 1,
\prod_{n\ge 1} N_n \le
N\Biggr\}.
\end{equation}

Then, let
$m\ge 2x_0+ 1$ be a temporarily fixed integer. We set, for $N\ge 1$,
\[
N_n = \biggl[\frac{a_n^pN^{1/m}}{(\prod_{1\le k\le m} a_k)^{
p/m}}\biggr],\qquad 1\le n\le m,\ N_n =1,\  n\ge m+1.
\]
The sequence $N_n,1\le n\le m$, is nonincreasing. This will ensure that
\[
   N_n\ge  1,\qquad 1\le n\le m.
\]
We wish to choose $m$ as a function of $N$ so that
\[
a_mN^{{1}/({pm})}\ge \Biggl(\prod_{1\le k\le m}a_k\Biggr)^{1/m}.
\]
Using $\log$-concavity, this is clearly satisfied provided that
%e4.22 ###
\begin{equation}\label{id1}
a_mN^{{1}/({pm})}\ge a_0\bigl((m+1)/2\bigr) =a\bigl((m+1)/2\bigr)
\end{equation}
[since  $(m+1)/2\ge x_0$]. Inequality~(\ref{id1}) becomes, by  taking logarithms,
%e4.23 ###
\begin{eqnarray}\label{id1bis}
&&\frac{m-1}{2} \log A +\frac{1}{p m}\log N
\vid
&&\qquad\ge   \frac{1}{\mu
p}\bigl(\log\bigl(\Gamma(m+1)\bigr) -\log\bigl(\Gamma\bigl(1+(m+1)/2\bigr)\bigr)\bigr).\nonumber
\end{eqnarray}

We will make use of the following classical inequality: for every
$t\ge 1/12$,
\[
0\le \log \bigl(\Gamma(t+1)\bigr)  -\log\bigl(\sqrt{2\pi}\bigr)- (t+1/2)\log t +t\le
1.
\]
Then, after some easy computations, one shows that
inequality~(\ref{id1bis}) is satisfied provided
\[
\frac{m-1}{2} \log A +\frac{1}{p  m}\log N \ge \frac{1}{\mu
p}\biggl(\frac {m}{2}\log m -\frac{m}{8}-\frac {1}{2} \log m +\frac {5}{ 2}\biggr).
\]
If one sets (this is probably  optimal)
\[
m=m(N):=\biggl\lceil2\sqrt{\mu\frac{\log N}{\log_2N}}\biggr\rceil,
\]
then the above inequality is satisfied,  as well as $m(N) \ge 2x_0+1$,
for every large enough $N$, provided that we increase the value of
$A$.   With
$N_n$ and $m$ settled as above and using the fact that
$\frac{x}{[x]}\le 2$ for every $x\ge 1$, we obtain
\[
\sum_{n\ge 1}\, \frac{a_n}{N^{1/p}_n}\le 2^{ 1/p} m
N^{-{1}/{(pm)}}\Biggl(\prod_{k=1}^ma_k\Biggr)^{1/m}
+\sum_{n\ge m+1}  a_n.
\]
On the one hand, $N_m\ge 1$ gives
\[
N^{-{1}/{(pm)}}\Biggl(\prod_{k=1}^ma_k\Biggr)^{1/m}\le  a_m.
\]
On the other hand, the $\log$-concavity and monotony   of the
function $a$ over $[x_0+1,\infty)$ (and the fact that $a'$ is nonzero)
imply that
\[
  \sum_{n\ge m+1}  a_n\le \bigg|\frac{a(x_0+1)}{a'(x_0+1)}\bigg|a_m = o(m a_m)
\]
(this follows from a straightforward adaptation of the proof of
Proposition~4.4 in {\cite{LUPA1}}, to which we refer the reader for
details). So we have
%e4.24 ###
\begin{eqnarray}\label{Coeur}
ma_m&=& m  \frac{A^m}{(m!)^{{1}/{(\mu p)}}} \nonumber\\
&\le &    \exp{\biggl(-\frac{1}{\mu p}m\log m +O(m)\biggr)}\\
&\le & C\exp{\biggl(-\frac{1}{p \mu}\sqrt{\mu \log
N\log_2N}\biggl(1+O\biggl(\frac{\log_3N}{\log_2N}\biggr)\biggr)\biggr)}.\nonumber
\end{eqnarray}

Note that $p\sqrt{\mu}= \sqrt{p\cdot p\mu}= \sqrt{(r+\delta)p}$. Finally,
this yields, in particular, that   for every $\varepsilon>0$,
\[
  \big\||K-\widehat K|_{L^p_{T}}\big\|_{r} =
O\biggl(\exp{\biggl(-\frac{1}{\sqrt{rp+\varepsilon}}\sqrt{\log
N\log_2N}\biggr)}\biggr).
\]
\end{step}
\begin{step}[(Compound case)]
Starting from (\ref{PoissonComp}),
it is natural to quantize
$(X_t)$ by setting
\[
\widehat{X}_t =\sum_{k=1}^{\widehat K_t} \widehat U_k,
\]
where $\widehat K$ is an $N^{(1)}$-quantization of the standard
Poisson  process $K$, as described in
Step~1, and, for every $n\ge 1$, $\widehat U_n$ is an $L^r$-optimal
$N_n^{(2)}$-quantization of $U_n$ with
$1\le N_1^{(2)}\times \cdots \times  N^{(2)}_n \cdots \le N^{(2)}$
and $N^{(1)}N^{(2)}\le N$. Then, setting
$\widehat K_t^U:=
\sum_{k=1}^{\widehat K_t}   U_k$ and $K_t^{\widehat U}:=
\sum_{k=1}^{ K_t}   \widehat U_k$, we obtain
\begin{eqnarray*}
|K^{\widehat U}- \widehat K^{\widehat U} |_{L^p_{T}}&\le &\sum_{n\ge
1} |\widehat U_k|\,\big|\mathbf{1}_{\{S_n\le
\lambda T\}}-\mathbf{1}_{\{\widehat S_n\le \lambda T\}}\big|_{L^p_{T}} \\
\end{eqnarray*}
so that
\begin{eqnarray*}
\big\||X-\widehat{K}^U|_{L^p_{T}}\big\|_{r}
&\le &\frac {1}{\lambda^{
1/p}}\sum_{n\ge 1} \|\widehat U_k\|_{r}\| S_n\wedge ( \lambda T)-
\widehat S_n\wedge (\lambda T)\|^{1/p}_{{r/p}} \\
&=&\frac{\sup_{n\ge 1}\|\widehat U_n\|_{r}}{\lambda^{1/p}}
\sum_{n\ge 1} \| S_n\wedge ( \lambda T)-
\widehat S_n\wedge (\lambda T)\|^{1/p}_{{r/p}},
\end{eqnarray*}
where  we have used the fact that the sequences $(U_n)$ and $(S_n)$ are
independent, as are $( \widehat U_n)$ and
$(S_n)$.  Using
\[
\| \widehat U_n\|_{r} \le \| U_n - \widehat U_n \|_{r}
+\|U_1\|_{r}=\big\| U_1 - \widehat U_1^{N^{(2)}_n} \big\|_{r} +\|U_1\|_{r}
\]
shows that $\sup_{n\ge 1}\|\widehat U_n\|_{r}<+\infty$. Hence, it
follows from Step~1 that, for every
$c<\frac{1}{\sqrt{pr}}$,
\[
\big\||X-
\widehat{K}^U|_{L^p_{T}}\big\|_{r}=O\bigl(\exp{\bigl(-c\sqrt{\log\bigl(N^{(1)}\bigr)\log_2\bigl(N^{(1)}\bigr)}\bigr)}\bigr).
\]
On the other hand, with obvious notation and using the fact that $(\widehat
U_n-U_n)$ and $(S_n)$ are independent, we have
\begin{eqnarray*}
\big\||X-K^{\widehat U}|_{L^p_{T}}\big\|_{r}&=& \big\||K^{U-\widehat
U}|_{L^p_{T}}\big\|_{r}\\
&\le& \sum_{n\ge 1} \|U_n-\widehat U_n\|_{r} \big\|\big|\mathbf{1}_{\{S_n\le \lambda\cdot\}}\big|_{L^p_{T}}\big\|_{r}\\
&=& \frac {1}{\lambda^{1/p}}\sum_{n\ge 1} \|U_n-\widehat
U_n\|_{r} \|(\lambda T-S_n)_+\|^{1/p}_{{r'}}\\
&\le& \frac {1}{\lambda^{1/p}}\sum_{n\ge 1} \|U_n-\widehat
U_n\|_{r} \frac{(\lambda T)^{1/p+n/r}}{(n!)^{1/r}}\\
&\le& C \sum_{n\ge 1} \|U_n-\widehat U_n\|_{r} \frac{(\lambda
T)^{n/r}}{(n!)^{1/r}}.
\end{eqnarray*}
It now follows from the (extended) Pierce lemma that
\begin{eqnarray*}
\big\||K^U-K^{\widehat U}|_{L^p_{T}}\big\|_{r}
&\le& C_{U_1,r} \sum_{n\ge 1}   \frac{(\lambda T)^{
n/r}}{(n!)^{1/r}N^{(2)}_n}\\
&=&O\biggl(\exp{\biggl(-\frac{1}{\sqrt{r}}\sqrt{\log\bigl(N^{(2)}\bigr)\log_2\bigl(N^{(2)}\bigr)}\biggr)}\biggr).
\end{eqnarray*}
The rate follows from  the resolution of the optimal bit allocation
problem~(\ref{bitalloc}) obtained by formally setting $\mu p =r$
and $p=1$.   Then note that, on the one hand,
\[
\big\||X-\widehat X|_{L^p_{T}}\big\|_{r}\le \big\||X-K^{\widehat
U}|_{L^p_{T}}\big\|_{r}+\big\||
K^{\widehat U}-\widehat K^{\widehat U}|_{L^p_{T}}\big\|_{r}
\]
and on the other hand
\[
\widehat K_t^{\widehat U}= \sum_{n\ge 1}
\widehat{U}_n^{N^{(2)}_n}\mathbf{1}_{\{\widehat S^{N_n^{(1)}}_n \le
\lambda t\}}
\]
can take at most
\[
\prod_{n\ge 1} N^{(1)}_nN^{(2)}_n \le N^{(1)}\times N^{(2)}\le N
\]
values. Let $c<\frac{1}{\sqrt{pr}}$. Setting $N^{(1)}=
[N^{{rc^2}/{(1+rc^2)}}]$, $N^{(2)}=
[N^{{1}/{(1+rc^2)}}]$ yields a rate
\[
\big\||X-\widehat X|_{L^p_{T}}\big\|_{r} =
O\biggl(\exp{\biggl(-\frac{1}{\sqrt{1/c^2+r}}\sqrt{\log(N)\log_2(N)}\biggr)}\biggr),
\]
that is,
\[
\forall\, \varepsilon>0\qquad \big\||X-\widehat X|_{L^p_{T}}\big\|_{r} =
O\biggl(\exp{\biggl(-\frac{1}{\sqrt{r(p+1+\varepsilon)}}\sqrt{\log(N)\log_2(N)}\biggr)}\biggr).
\]\hspace{337pt}\ \qed\noqed
\end{step}\noqed
\end{pf*}\noqed

%s4.3 ###
\subsection{\texorpdfstring{Proof of Theorem \textup{\protect\ref{Levyrate}}}{Proof of Theorem 2}}\label{finalstep}
Any L\'evy process $X$ can be decomposed as the sum $X= X^{(1)}+X^{(2)}$
of two
(independent) L\'evy processes, one having bounded jumps and the other
being  a compound Poisson process,
according to the decomposition of its L\'evy measure
%
%e4.25 ###
\begin{eqnarray}\label{nu12}
\quad&&\nu(dx)= \nu^{(1)}(dx)+\nu^{(2)}(dx)
\vid
&&\qquad\mbox{ with } \nu^{(1)}(dx):=
\mathbf{1}_{\{|x|\le 1\}}\nu(dx)
\mbox{ and }
\nu^{(2)}(dx):=\mathbf{1}_{\{|x|> 1\}}\nu(dx).\nonumber
\end{eqnarray}
Assume that $r^*> 1$. It is then clear that, for every $r, p\!\in(0,
r^*)$,
%e4.26 ###
\begin{eqnarray}\label{identSum1}
e_{N,r}(X, L^p_{T}) &\le& C_{r,p,T} e_{[\sqrt{N}]^2,r'}(X,
L^{r'}_{T})
\vid
&\le& C_{r,p,T}\bigl(e_{[\sqrt{N}],r'}\bigl(X^{(1)},
L^{r'}_{T}\bigr)+e_{[\sqrt{N}],r'}\bigl(X^{(2)}, L^{r'}_{T}\bigr)\bigr),\nonumber
\end{eqnarray}
where $r'\!=\!r \vee p\vee 1$. It now follows from
Proposition~\ref{Poisson} that $e_{N,r'}(X^{(2)},
L^{r'}_{T})\!=\!o(e_{N,r'}(X^{(1)},
L^{r'}_{T}))$ so that
\[
e_{N,r}(X, L^p_{T}) \le C'_{r,p,T }e_{[\sqrt{N}],r'}\bigl(X^{(1)}, L^{r'}_{T}\bigr).
\]
Now, using the fact that $\underline\ell$ has slow variations at $0$, we can derive that
\[
e_{[\sqrt{N}],r'}\bigl(X^{(1)}, L^{r'}_{T}\bigr)=O\bigl(e_{N,r'}\bigl(X^{(1)}, L^{r'}_{T}\bigr)\bigr).
\]
Proposition~\ref{sous-theta}  completes the proof of
Theorem~\ref{Levyrate}. When $r^*\le 1$,  we use
\begin{eqnarray*}
e_{N,r}(X, L^p_{T})^r &\le& C'_{r,p,T} e_{[\sqrt{N}]^2,r'}(X,
L^{r'}_{T})\\
&\le& C'_{r,p,T} \bigl(e_{[\sqrt{N}],r'}\bigl(X^{(1)},
L^{r'}_{T}\bigr)^{r'} + e_{[\sqrt{N}],r'}\bigl(X^{(2)},
L^{r'}_{T}\bigr)^{r'}\bigr),
\end{eqnarray*}
with $r'=r\vee p<1$ (based on the pseudo-triangular inequality
satisfied by\break $L^s$-pseudo-norms when $s<1$).

%s5 ###
\section{Further  results for L\'evy processes}\label{toto}
%s5.1 ###
\subsection{An exact  rate for L\'evy processes with a Brownian component}\label{LevyB}
In that case, the quantization rate of the Brownian motion controls the
global rate of convergence.

\begin{proposition}\label{LevywithB} Let $X$ be a L\'evy process with  a
nonvanishing Brownian component. Let  $r^*=r^*(X)$, defined
by~(\ref{thetahaut}). Then
\[
\forall\,r,\, p\in (0,r^*\wedge 2)\qquad e_{N,r}(X,L^p_{T})=
O((\log N)^{-1/2})
\]
and
\[
\forall\,r,\, p\in (0,+\infty)\qquad e_{N,r}(X,L^p_{T})=
\Omega(e_{{N,r}}(W, L^{p}_{T})).
\]
In particular, $\forall\,r\!\in(0,+\infty),\,\forall\, p\!\in
[1,+\infty)$, $ e_{N,r}(X,L^p_{T})= \Omega((\log N)^{-1/2})$.
\end{proposition}

\begin{pf}
We can decompose $X$ as
\[
X=cW+X^{(1)}+X^{(2)},
\]
where $X^{(i)}$, $i=1,2$, have $\nu^{(i)}$ as L\'evy measure, as
defined in~(\ref{nu12}) in the above proof of Theorem~\ref{Levyrate}.
Then, if $r^*>1$ and $r,\,p\!\in(0,r^*\wedge 2)$, we can easily check
that, for every $N\ge 1$,
\begin{eqnarray*}
e_{{N,r'}}(X,L^p_{T}) &\le &e_{{[\sqrt[3]{N}]^3,r'}}(X,L^{r'}_{T}) \\
&\le&
e_{{[\sqrt[3]{N}],r'}}(c\,W,L^{r'}_{T})+e_{{[\sqrt[3]{N}],r'}}\bigl(X^{(1)},L^{r'}_{T}\bigr)
+ e_{{[\sqrt[3]{N}],r'}}\bigl(X^{(2)},L^{r'}_{T}\bigr),
\end{eqnarray*}
where $r'\!=\!r\vee p\vee  1$.  It follows from
Proposition~\ref{Poisson} (see the remark immediately  below) that
$e_{N,r'}(X^{(2)}, L^{r'}_{T})\!=\!o(e_{N,r'}(W, L^{r'}_{T}))$.
Now, $\int_{\mathbb{R}\setminus \{0\}}x^2\nu^{(1)}(dx)<+\infty$, hence, by
Millar's
lemma,
\[
\mathbb{E} \sup_{s\in [0,t]}\big|X^{(1)}_s\big|^2 \le Ct.
\]
We can then easily derive from~(\ref{MillarRate}) (or directly from
Theorem~\ref{Upper1}) that\break $e_{{N,r'}}(X^{(1)}, L^{r'}_{T}) =
O((\log N)^{-1/2})$. This yields the announced upper bound since
$e_{{N,r'}}(W, L^{r'}_{T}) = O((\log N)^{-1/2})$. If $r^*\le
1$, we proceed as above, using the pseudo-triangular inequality  for
$L^s$-pseudo-norms (with $ r '\! = \! r \vee  p < 1 $).

As concerns the lower bound, note that if $Y$ and $Z$ are
$L^r_{T}$-valued independent random vectors, then for every $r,\,p>0$,
\begin{eqnarray*}
\bigl( e_{{N,r}}(Y+Z, L^{p}_{T}) \bigr)^r  &  =  &
\inf_{\alpha \subset L^{p}_{T},\,\operatorname{card}(\alpha)\le N} \int
\mathbb{E}\min_{a\in \alpha}|Y-z-a|_{L^{p}_{T}} ^r \mathbb{P}_{Z}(dz)\\
&\ge & \int_{L^{p}_{T}}  \inf_{\alpha \subset L^{p}_{T}, \,
\operatorname{card}(\alpha)\le N} \mathbb{E}\min_{a\in \alpha}|Y-z-a|_{L^{p}_{T}} ^r
\mathbb{P}_{Z}(dz)\\
&=& (e_{{N,r}}(Y,L^{p}_{T}))^r
\end{eqnarray*}
so that
\[
e_{{N,r}}(Y+Z, L^{p}_{T}) \ge \max(e_{{N,r}}(Y, L^{p}_{T}) ,
e_{{N,r}}(Z, L^{p}_{T})).
\]
This holds true, by induction, for any finite sum of independent random
variables. In particular,
\[
e_{{N,r}}(X,L^p_{T})\ge  e_{{N,r}}(cW,
L^{p}_{T})=ce_{{N,r}}(W, L^{p}_{T}).
\]
This  completes the proof.\quad\qed
\end{pf}

%s5.2 ###
\subsection{Subordinated L\'evy processes}\label{suboLevy}
We  now consider subordination of the Brownian motion, that is, L\'evy
processes of the form
\[
X_t = W_{A_t}, \qquad t\ge 0,
\]
where $W$ denotes a standard Brownian motion  and $A$ a subordinator
independent of $W$. A subordinator is a nondecreasing
(hence nonnegative) L\'evy process. What follows is borrowed from~\cite{BER0}. Its
L\'evy--Khintchine characteristics $(a,\sigma^2,\nu_{{\!A}})$ satisfy
$\sigma^2=0$, $\nu_{{\!A}}((-\infty,0))=0$,
$\int_0^1x\nu_{{\!A}}(dx)<+\infty$ and
$\gamma:=a-\int_0^1x\nu_{{\!A}}(dx)\ge 0$ [so that $\underline
\theta(A)\le 1]$. Consequently, a subordinator is of the form
\[
A_t = \gamma t +\sum_{s\le t} \Delta A_s, \qquad t\ge 0.
\]
Its Laplace transform is given by $\mathbb{E}\, e^{-uA_t} = e^{-t\Phi(u)}$
with, for every $u\ge 0$,
%e5.1 ###
\begin{eqnarray}\label{Phisub}
\Phi(u)&=&\gamma u+\int_0^{+\infty}
(1-e^{-ux})\nu_{{A}}(dx)
\vid
&=&
\gamma u +u\int_0^{+\infty} e^{-ux}\underline \nu_{{A}}(x)\,dx,\nonumber
\end{eqnarray}
where $\underline \nu_{{A}}(x)= \nu_{{A}}((x,+\infty))$ denotes
the tail of the L\'evy measure $\nu_{{A}}$ and $
\lim_{u\to+ \infty}\frac{\Phi(u)}{u} = \gamma$.
Furthermore, for every $t\ge 0$ and $u\in \mathbb{R}$,
\[
\mathbb{E}(e^{iuX_t})= \mathbb{E}\bigl(e^{(-{u^2}/{2})A_t}\bigr) =
\exp{\biggl(-\frac{u^2}{2}\gamma-\int_0^{+\infty}\bigl(1-e^{-({u^2}/{2})x}\bigr)\nu_{{A}}(dx)
\biggr)}
\]
so that we can easily derive that
\[
\nu_{{X}}(f) = \int_{(0,\infty)}\mathbb{E}\bigl(f\bigl(\sqrt{x}\,Z\bigr)\bigr)\nu_{{A}}(dx)
\]
[with $Z\sim\mathcal{N}(0;1)$]  and that $X$ has a Brownian component if
and only if $\gamma>0$ (see also~\cite{SAT}, page 198).

The small deviation of subordinator has been extensively investigated
in~\cite{LISH}. It is there established that if $
\liminf_{u\to +\infty} \frac{\Phi(u)}{\log u} >0$, then
%e5.2 ###
\begin{equation}\label{sblbsub}
\qquad\forall\, p\!\in [1,+\infty) \qquad -\log\bigl(\mathbb{P}\bigl(|X|_{L^p_{T}}\le
\varepsilon\bigr)\bigr)\approx \Phi(\varepsilon^{-2})\qquad\mbox{as }
\varepsilon \to 0.
\end{equation}
These processes preserve a Gaussian feature which will be the key to
estimate their quantization rate: they satisfy the Anderson
inequality, as briefly recalled in the lemma below.
\begin{lema}\label{Andersbis}
A subordinated L\'evy process is unimodal for every $L^p_{T}$-norm,
for every $p\!\in[1,+\infty)$. The result still holds if one replaces
$W$ by, for example, any pathwise continuous centered Gaussian process
(e.g., fractional Brownian motion, etc.).
\end{lema}
\begin{pf}
Using the fact that $A$ and $W$ are independent, it
suffices to show that for every nondecreasing function $a\dvtx[0,T]\to
[0, \alpha(T)]$, $a(0)=0$, and every $x\!\in L^p_{T}$,
\[
\mathbb{P}\biggl(\int_0^T \big| W_{a(s)}-x(s)\big|^p\,ds \le \varepsilon^p\biggr)\le
\mathbb{P}\biggl(\int_0^T \big| W_{a(s)}\big|^p\,ds \le
\varepsilon^p\biggr),\qquad\varepsilon>0.
\]
It is clear that $(W_{a(t)})_{t\in[0,T]}$ is a centered
(bimeasurable) Gaussian  process and has sample paths in $L^p_{T}$
a.s. Hence,  $(W_{a(t)})_{t\in[0,T]}$ can be seen as an
$L^p_{T}$-valued centered Gaussian random vector and the assertion
follows from the Anderson inequality.%\rightqed
\end{pf}

We now make the connection between Blumenthal--Getoor indices of $X$
and $A$ (and between the finiteness of  moments).

\begin{lema} \label{BGsublevy}
$\underline \theta(X)=2\,\underline
\theta(A)$ and $r^*(X)=2\,r^*(A)$.
\end{lema}

\begin{pf}
As a consequence of the expression for $\nu_{{\!X}}$, we check
that for every $\theta\!\in (0,2]$,
\begin{eqnarray*}
\int_{\{|x|\le 1\}}|x|^\theta \nu_{X}(dx)&=& \int
u^{{\theta}/{2}}\int_{\{|y|\le 1/\sqrt{u}\}}|y|^\theta
e^{-{y^2}/{2}}\frac{dy}{\sqrt{2\pi}}\\
&=&  \int_{\{u>0\}}
u^{{\theta}/{2}} \nu_{{A}}\ell_\theta(u)(du),
\end{eqnarray*}
where $\ell_\theta(u)> 0$ when $u\!>\!0$
and   $\lim_{u\to 0}\ell_\theta(u) \!=\! C_\theta\!\in (0,+\infty)$.
Hence, the first equality follows. As concerns the second equality, $r^*(X)$
coincides with the (absolute) moments of $X$, so it is obvious that
\[
\mathbb{E}(|X_t|^r)=\mathbb{E}(|W_{{A_t}}|^r)= \mathbb{E}(A^{2r}_t).
\]
Consequently,  $\mathbb{E}(|X_t|^r)<+\infty$ iff $\mathbb{E}(A^{2r}_t)<+\infty$ so
that  $r^*(X)=2r^*(A)$.\quad\qed
\end{pf}

As concerns upper bounds, we cannot apply Theorem~\ref{Levyrate}
since a subordinated L\'evy process may have  a Brownian component.
Therefore, we must return Theorem~\ref{Upper1}.
\begin{proposition}
\textup{(a)} If $\gamma>0$, then
\[
\forall\, r,\,p\in\bigl(0,r^*(X)\wedge 2\bigr))\qquad e_{{N,r}}(X,
L^p_{T}) = O((\log N)^{-1/2}).
\]
\textup{(b)} If  $\underline \theta(A)\in(0,1)$, $\gamma=0$ and
$\nu_{{\!A}}(dx)\mathbf{1}_{\{0<x\le \eta\}}\le c\mathbf{1}_{\{0<x\le \eta\}}\frac{dx}{x^{1+\underline \theta(A)}}$ for some
real constants $c,\eta>0$, then
\[
\forall\, r,\,p\!\in(0,\underline \theta(X) \wedge r^*(X))\qquad
e_{{N,r}}(X, L^p_{T}) = O\bigl((\log N)^{-{1}/{(\underline
\theta(X))}}\bigr).
\]
\end{proposition}

\begin{pf}
(a)   follows from
Proposition~\ref{LevywithB} since $X$ has a Brownian component.

(b) Let $\rho < 2\,(\underline \theta(A)\wedge r^*(A))$.
First, note that $\mathbb{E}(|X_t|^\rho)=  \mathbb{E} A_t^{\rho/ 2} \le C
t^{{\rho}/{2\underline \theta(A)}}$, $\rho\le 2 (\underline
\theta(A)\wedge r^*(A))=\underline \theta(X)\wedge r^*(X)$ (by
Lemma~\ref{Quantifrate1}
applied to $A$). The result then follows from Theorem~\ref{Upper1}.\qed
\end{pf}

The following lower bounds follow from Lemma~\ref{Andersbis} and
inequality~(\ref{sblowerbound}) (see the remark immediately after
Theorem~\ref{Upper1}). The main point to be noted is that \textit{the
upper and lower bounds obtained match, providing an exact
quantization rate for subordinated L\'evy processes}.
\begin{proposition}
\textup{(a)} If $\gamma>0$, then
\[
\forall\, r\!\in(0,+\infty),\; \forall\, p\!\in [1,+\infty)\qquad
e_{{N,r}}(X, L^p_{T}) = \Omega((\log N)^{-1/2}).
\]

\textup{(b)} If $\gamma=0$, $\underline \theta(A)>0$  and $\mathbf{1}_{\{0<x\le \eta\}}\nu_{{\!A}}(dx)\ge c\mathbf{1}_{\{0<x\le
\eta\}}\frac{dx}{x^{1+\underline \theta(A)}}$ for some real constants
$c,\, \eta>0$, then
\[
\forall\, r\!\in(0,+\infty),\; \forall\, p\!\in [1,+\infty)\qquad
e_{{N,r}}(X, L^p_{T}) = \Omega\bigl((\log N)^{{1}/{\underline
\theta(X)}}\bigr).
\]
\end{proposition}
\begin{pf}
(a)  follows from
Proposition~\ref{LevywithB} since $X$ has a Brownian component.

(b) It follows from the assumption made on $\nu_{{A}}$
that $\underline \nu_{{A}}(x) \ge c\int_x^\eta \xi^{-\underline
\theta-1}\,d\xi \ge\break ~\kappa
x^{-\underline \theta}$ for
$x\in (0, \eta/2]$. Hence, it follows from~(\ref{Phisub}) that
\[
\Phi(u) \ge cu \int_0^{\eta/2} e^{-ux}\underline \nu_{{\!A}}(x)\,dx
=cu^{\underline \theta} \int_0^{ u\eta/2} e^{-y}y^{-\underline
\theta}\,dy\ge c' u^{\underline \theta}
\]
for  large enough $u$ (with an appropriate real constant $c'>0$). We
conclude by combining~(\ref{sblowerbound}) and~(\ref{sblbsub}) since
$X$ is strongly unimodal.\qed
\end{pf}

\begin{examples*}
If $A$ is a tempered
$\alpha$-stable process with L\'evy measure, then
\[
\nu_{{\!A}}(dx)=
\frac{2^\alpha\alpha}{\Gamma(1-\alpha)}x^{-(\alpha+1)}\exp{\biggl(-\frac{1}{2}\delta^{
1/\alpha}\biggr)}\mathbf{1}_{(0,\infty)}(x)\, dx,
\]
with $\alpha\!\in (0,1)$, $\delta>0$, $\gamma=0$, so that $\underline
\theta(A)=\alpha$ and $r^*(A) +\infty$. We the obtain
\[
\forall\, r\!\in(0,2\alpha),\; \forall\,p\!\in [1,2\alpha)\qquad
e_{{N,r}}(X,L^p_{T})\approx (\log N)^{-{1}/{(2\alpha)}}.
\]

Assume that $\underline \theta(A)\in (0,1)$ and that
the function $\Phi$ is regularly varying at $\infty$ with index
$\alpha \!\in (0,1)$ such that
\[
\Phi(x)\sim c x^{\alpha}(\log(x))^c\qquad \mbox{as } x\to
\infty,
\]
for some real constant $c>0$. Since $\alpha<1$, we have $\gamma=0$. Then
\[
\Gamma(1-\alpha)\underline \nu(x) \sim \Phi(1/x)\qquad \mbox{ as }x\to 0
\]
(see~\cite{BER0}) so that $\underline \nu$ is regularly varying at
zero with index $-\alpha$. By Theorem~\ref{Levyrate}, $\underline
\theta(A)=\alpha$. Set
  \[
  \Psi(x) = x^{{1}/{(2\underline \theta(A))}} (\log
x)^{-{c}/{(2\underline \theta(A))}}
  \]
for large enough $x>0$. Then $\Psi\circ \Phi(x) \sim c\sqrt{x}$
as $x\to \infty$ so that $\Psi\circ \Phi(1/\varepsilon^2)\sim c
\varepsilon^{-1}$ as $\varepsilon \to 0$. Thus,
\begin{eqnarray*}
&&\forall\,r>0, \forall\,p\!\in [1,+\infty)\\
&&\qquad e_{{_N},r}(X,L^p_{T})=\Omega\bigl((\log N)^{-{1}/{(2\underline
\theta(A))}}(\log\log N)^{-{c}/{(2\underline \theta(A))}}
\bigr).
\end{eqnarray*}
On the other hand, by Lemma~\ref{Quantifrate1} and remark below
Theorem~\ref{Upper1}, in the case $c>0$,
\[
  \mathbb{E}\, A_t^{{\rho}/{2}} \le C t^{{{\rho}/({2\underline
\theta(A)})}}(-\log t)^c,\qquad \rho/2< \underline \theta(A)\wedge
r^*(A)
\]
so that
\begin{eqnarray*}
&&\forall\, r,\, p\in \bigl(0,\underline\theta(X)\wedge r^*(X)\bigr),\\
&&\qquad
e_{{N,r}}(X,L^p_{T})=
O\bigl((\log N)^{-{1}/{\underline \theta(X)}}(\log\log N)^{{c}/{\rho}}
\bigr),\qquad \rho< \underline \theta(X)\wedge r^*(X).
\end{eqnarray*}
In the case $r^*(X)\ge \underline \theta(X)$, this matches the lower
bound up to a $O(\log\log N)^{\varepsilon}$ term, $\varepsilon>0$.
\end{examples*}

\section*{Acknowledgments} We thank Sylvain Delattre, G\'erard
Kerkyacharian and
Dominique Picard for helpful discussions.

\printaddresses

\begin{thebibliography}{99}
%b1 ###
%
%b2 ###
\bibitem{BER0} \textsc{Bertoin, J.} (1999).
Subordinators: Examples and applications.  \textit{Lectures on Probability
Theory and Statistics \textup{(}Saint-Flour\textup{,} 1997\textup{)}}.  \textit{Lecture Notes in
Math.} \textbf{1717} 1--91. Springer, Berlin.
\MR{1746300}
%
%b3 ###
\bibitem{BER} \textsc{Bertoin, J.} (1996). \textit{L\'evy Processes.}
Cambridge Univ. Press.
\MR{1406564}
%
%b4 ###
\bibitem{BIGOTE}\textsc{Bingham, N. H., Goldie, C. M.} and  \textsc{Teugels, J. L.} (1987).
\textit{Regular Variation}. Cambridge Univ. Press.
\MR{0898871}
%
%b5 ###
\bibitem{BOLE}\textsc{Bouleau, N.}  and \textsc{L\'epingle, D.} (1994).
\textit{Numerical Methods for Stochastic Processes}. Wiley, New York.
\MR{1274043}
%
%b6 ###
\bibitem{BLGE} \textsc{Blumenthal, R. M.} and \textsc{Getoor, R. K.} (1961). Sample
functions of stochastic processes with stationary independent
increments.
\textit{J. Math. Mech.} \textbf{10}  493--516.
\MR{0123362}
%
%b7 ###
\bibitem{DEFEMASC} \textsc{Dereich, S., Fehringer, F., Matoussi, A.}
and
\textsc{Scheutzow, M.} (2003).  On the link
between small ball probabilities and the quantization problem for
Gaussian measures on Banach spaces. \textit{J. Theoret. Probab.} \textbf{16} 249--265.
\MR{1956830}
%
%b8 ###
\bibitem{Dereich} \textsc{Dereich, S.} (2005). The coding complexity of
diffusion processes under $L^p([0,1])$-norm distortion. \textit{Stochastic
Process.
Appl.} To appear.
%
%b9 ###
\bibitem{Dereichetal} \textsc{Dereich, S.} and \textsc{Scheutzow, M.} (2006). High
resolution quantization and entropy coding for fractional Brownian
motions. \textit{Electron. J. Probab.}  \textbf{11} 700--722.
\MR{2242661}
%
%b10 ###
\bibitem{EBKE}  \textsc{Eberlein, E.} and \textsc{Keller, U.} (1995). Hyperbolic
distributions in finance. \textit{Bernoulli} \textbf{1} 281--299.
%
%b11 ###
\bibitem{EMMA} \textsc{Embrechts, P.} and \textsc{Maejima, M.} (2002).
\textit{Self-Similar Processes}. Princeton Univ. Press. \MR{1920153}
%
%b12 ###
\bibitem{IEEE} \textsc{Gersho, A.} and  \textsc{Gray, R. M.} (1982). Special issue on
quantization. \textit{IEEE Trans. Inform.
Theory} \textbf{29}.
%
%b13 ###
\bibitem{GEGR} \textsc{Gersho,  A.} and \textsc{Gray, R. M.} (1992). \textit{Vector
Quantization and Signal Compression}. Kluwer, Boston.
%
%b14 ###
\bibitem{GRLU}
\textsc{Graf, S.} and \textsc{Luschgy, H.} (2000). \textit{Foundations of Quantization
for Probability Distributions}. Springer, Berlin.
\MR{1764176}
%%
%b15 ###
\bibitem{GRLUPA1} \textsc{Graf, S., Luschgy, H.} and \textsc{Pag\`es, G.} (2003).
Functional quantization and small ball probabilities for Gaussian
processes.  \textit{J. Theoret. Probab.}  \textbf{16} 1047--1062.
\MR{2033197}
%
%b16 ###
\bibitem{GRAF3} \textsc{Graf,  S., Luschgy,  H.} and \textsc{Pag\`es,  G.} (2007).
Optimal quantizers for Radon random vectors in a Banach space.
  \textit{J. Approx. Theory} \textbf{144} 27--53.
\MR{2287375}
%
%b17 ###
\bibitem{JASH} \textsc{Jacod, J.} and \textsc{Shiryaev, A.} (2003). \textit{Limit
Theorems for Stochastic Processes}, 2nd ed. Springer,
Berlin.
\MR{1943877}
%
%b18 ###
%for fractional stable processes. \textit{Ann. Inst. H. Poincar\'e
%Probab. Statist. } \textbf{41} 725--752.
%
%b19 ###
\bibitem{LISH} \textsc{Linde, W.} and \textsc{Shi, Z.} (2004). Evaluating the small
deviation probabilities for subordinated L\'evy
processes.  \textit{Stochastic Process. Appl.}  \textbf{113} 273--287.
\MR{2087961}
%
%b20 ###
\bibitem{LUPA1}
\textsc{Luschgy, H.} and \textsc{Pag\`es, G.} (2002). Functional quantization of
Gaussian processes. \textit{J. Funct. Anal.}  \textbf{196} 486--531.
\MR{1943099}
%
%b21 ###
\bibitem{LUPA2}
\textsc{Luschgy, H.} and  \textsc{Pag\`es, G.} (2004). Sharp asymptotics of the
functional quantization problem
for Gaussian processes. \textit{Ann. Probab.} \textbf{32}
1574--1599.
\MR{2060310}
%
%b22 ###
\bibitem{LUPA4}
\textsc{Luschgy, H.} and  \textsc{Pag\`es, G.} (2006). Functional quantization of a
class of Brownian diffusions: A constructive approach.  \textit{Stoch.
Process.
Appl.} \textbf{116} 310--336.
\MR{2197980}
%
%b23 ###
\bibitem{MI} \textsc{Millar, P. W.} (1971).   Path behaviour of processes
with stationary independent increments.  \textit{Z.
Wahrsch. Verw. Gebiete} \textbf{17} 53--73.
\MR{0324781}
%
%b24 ###
\bibitem{PAPR2}
\textsc{Pag\`es, G.} and \textsc{Printems, J.} (2005).   Functional quantization for
numerics with an application to
option pricing.  \textit{Monte Carlo Methods
Appl.} \textbf{11} 407--446.
\MR{2186817}
%
%b25 ###
\bibitem{POST} \textsc{P\"otzelberger, K.} and \textsc{Strasser, H.} (2001). Clustering
and quantization by MSP-partitions.  \textit{Statist. Decisions}   \textbf{19} 331--371.
\MR{1884124}
%
%b26 ###
\bibitem{SAT} \textsc{Sato, K. I.} (1999). \textit{L\'evy Processes and
Infinitely Divisible Distributions}. Cambridge Univ. Press.
\MR{1739520}
%
%b27 ###
\bibitem{SATA} \textsc{Samorodnitsky, G.} and \textsc{Taqqu, M. S.} (1994). \textit{Stable
non-Gaussian Random Processes}. Chapman and Hall, New York.
\MR{1280932}
%
%b28 ###
\bibitem{SI} \textsc{Singer, I.} (1970).\textit{ Basis in Banach Spaces.} \textit{I}.
Springer, New
York.
\MR{0298399}
%
%b29 ###
\bibitem{TAKI} \textsc{Tarpey, T.} and \textsc{Kinateder, K. K. J.} (2003). Clustering
functional data. \textit{J. Classification} \textbf{20} 93--114.
\MR{1983123}
%
%b30 ###
\bibitem{TAPEOG} \textsc{Tarpey, T., Petkova, E.} and \textsc{Ogden, R. T.} (2003).
Profiling placebo responders by self-consistent partitioning of
functional
data. \textit{J. Amer. Statist. Assoc.} \textbf{98} 850--858.
\MR{2055493}
%
\bibitem{WI} \textsc{Wilbertz, B.} (2005). Computational aspects of
functional quantization for Gaussian measures and applications.
Diploma thesis, Univ. Trier.
\end{thebibliography}
\end{document}